\documentclass[10pt]{amsart}
\usepackage[sort,square,compress,comma, numbers]{natbib}
\allowdisplaybreaks[4]
\usepackage{amssymb,color}
\usepackage[colorlinks=true, citecolor=blue, linkcolor=blue]{hyperref}
\usepackage{amsthm}

\definecolor{c20}{rgb}{0.,0.7,0.}
\definecolor{c30}{rgb}{0.,0.,1.}
\definecolor{c40}{rgb}{1,0.1,0.7}
\definecolor{c50}{rgb}{1,0,0}
\definecolor{c60}{rgb}{0,0.9,0.1}

\newcommand{\ve}{\varepsilon}

\newcommand{\abs}[1]{\lvert #1 \rvert}
\newcommand{\ABs}[1]{ \biggl \lvert #1 \biggr \rvert}
\newcommand{\norm}[1]{\lVert #1 \rVert}

\newcommand{\E}[1]{\mathbb{E}\left\{ #1 \right\}}
\newcommand{\EE}[1]{\mathbb{E}\left\{ #1 \right\}}

\newcommand{\pk}[1]{\mathbb{P} \left\{ #1 \right\} }

\newcommand{\R}{\mathbb{R}}

\newcommand{\N}{\mathbb{N}}
\newcommand{\inr}{\in \R}
\newcommand{\inn}{\in \N}

\newcommand{\ldot}{,\ldots,}

\newcommand{\limit}[1]{\lim_{#1 \to   \infty}}

\def\Ha{ \mathcal{H}_{\alpha}}

\newcommand{\BQN}{\begin{eqnarray}}
\newcommand{\EQN}{\end{eqnarray}}
\newcommand{\BQNY}{\begin{eqnarray*}}
\newcommand{\EQNY}{\end{eqnarray*}}

\def\K#1{\textcolor{magenta}{#1}}
\def\K#1{{#1}}

\def\bqny#1{{ \begin{eqnarray*} #1 \end{eqnarray*}}}
\def\bqn#1{ { \begin{eqnarray} #1 \end{eqnarray}}}

\newcommand{\BS}{\begin{sat}}
\newcommand{\ES}{\end{sat}}
\newcommand{\BT}{\begin{theo}}
\newcommand{\ET}{\end{theo}}
\newcommand{\BK}{\begin{korr}}
\newcommand{\EK}{\end{korr}}

\newcommand{\BD}{\begin{de}}
\newcommand{\ED}{\end{de}}

\newcommand{\BIT}{\begin{itemize}}
\newcommand{\EIT}{\end{itemize}}
\newcommand{\BDI}{\begin{description}}
\newcommand{\EDI}{\end{description}}

\newcommand{\BRM}{\begin{remark}}
\newcommand{\ERM}{\end{remark}}

\newcommand{\BEL}{\begin{lem}}
\newcommand{\EEL}{\end{lem}}

\newtheorem{theo}{Theorem}[section]
\newtheorem{sat}[theo]{Proposition}
\newtheorem{de}[theo]{Definition}
\newtheorem{lem}[theo]{Lemma}

\newtheorem{korr}[theo]{Corollary}
\newtheorem{remark}[theo]{Remark}

\newcommand{\netheo}[1]{{Theorem \ref{#1}}}

\newcommand{\prooftheo}[1]{ \textbf{Proof of Theorem} \ref{#1} }
\newcommand{\proofprop}[1]{\textbf{Proof of Proposition} \ref{#1}}

\newcommand{\proofkorr}[1]{\textbf{Proof of Corollary} \ref{#1}}

\newcommand{\proofrem}[1]{\textbf{Proof of Remark} \ref{#1}}
\newcommand{\COM}[1]{}

\newcommand{\QED}{\hfill $\Box$ \\}

\def\rw{\rightarrow}

\def\IF{\infty}

\def\LT{\left}
\def\RT{\right}

\def\eHH#1{\textcolor{c50}{#1}}

\def\eHHH#1{\textcolor{c20}{#1}}
\def\eHH#1{#1}
\def\eHHH#1{#1}

\def\ehd#1{\textcolor{black}{#1}}
\def\gE#1{\textcolor{c50}{#1}}

\def\gE#1{{#1}}
\def\eM#1{\textcolor{c50}{#1}}
\def\gE#1{{#1}}
\def\eM#1{\textcolor{c50}{#1}}
\def\HEH#1{{#1}}

\def\eM#1{{#1}}

\def\KD#1{\textcolor{black}{#1}}
\def\KD#1{{#1}}
\def\cl#1{\textcolor{black}{#1}}

\def\k#1{\textcolor{black}{#1}}
\def\rd#1{\textcolor{black}{#1}}
\def\rdd#1{\textcolor{black}{#1}}
\def\kd#1{\textcolor{black}{#1}}
\def\kk#1{\textcolor{black}{#1}}
\def\k#1{{#1}}
\def\rrd#1{\textcolor{black}{#1}}
\def\vf{\sigma^2}

\def\gtu{g_{u,\tau_u}}
\def\xitu{\xi_{u,\tau_u}}

\def\chitu{\chi_{u,\tau_u}}
\def\zetatu{\zeta_{u,\tau_u}}
\def\vtu{\sigma_{u,\tau_u}}
\def\bD{E}
\def\sigxiu{\vtu}

\def\coe{C_0(E)}

\begin{document}

\title[Tails of homogenous functionals of Gaussian  fields ]
{Uniform Tail Approximation of homogenous functionals of Gaussian fields}

\author{Krzysztof D\c{e}bicki}
\address{Krzysztof D\c{e}bicki, Mathematical Institute, University of Wroc\l aw, pl. Grunwaldzki 2/4, 50-384 Wroc\l aw, Poland}
\email{Krzysztof.Debicki@math.uni.wroc.pl}
\author{Enkelejd Hashorva}
\address{Enkelejd Hashorva, Department of Actuarial Science, University of Lausanne, UNIL-Dorigny 1015 Lausanne, Switzerland}
\email{enkelejd.hashorva@unil.ch}
\author{Peng Liu}
\address{Peng Liu, Department of Actuarial Science, University of Lausanne, UNIL-Dorigny 1015 Lausanne, Switzerland and Mathematical Institute, University of Wroc\l aw, pl. Grunwaldzki 2/4, 50-384 Wroc\l aw, Poland}
\email{peng.liu@unil.ch}

\bigskip

\date{\today}
 \maketitle

\bigskip
{\bf Abstract:} Let $X(t),t\inr^d$ be a centered Gaussian random field with continuous trajectories and set $\xi_u(t)= X(f(u)t),t\inr^d$ with $f$ some positive function. Classical results establish the tail
asymptotics of $\pk{ \Gamma(\xi_u) > u}$ as $u\to \IF$
with $\Gamma(\xi_u)= \sup_{t \in [0,T ]^d} \xi_u(t),T>0$ by requiring that  $f(u)$ tends to 0 as $u\to \IF$ with speed controlled by the local behaviour of the correlation function of $X$.
Recent research shows that for applications more general 
functionals than supremum should be considered and the Gaussian  field can depend also on some additional parameter $\tau_u \in K$, say $\xi_{u,\tau_u}(t),t\inr^d$. In this contribution we derive uniform approximations of $\pk{ \Gamma(\xi_{u,\tau_u})> u}$ with respect to $\tau_u$ in some index set $K_u$, as $u\to\infty$.  Our main result have important theoretical implications; two applications are already included in \cite{KEP2015, KEP2016}. In this paper we present three additional ones, namely i) we derive uniform upper bounds for the probability of double-maxima, ii) we extend Piterbarg-Prisyazhnyuk theorem to some large classes of homogeneous functionals of centered Gaussian fields $\xi_{u}$, and iii) we show the finiteness of generalized Piterbarg constants.

{\bf Key Words}: fractional Brownian motion; supremum of Gaussian random fields; stationary processes; double maxima; uniform double-sum method; generalized Piterbarg constants.

{\bf AMS Classification:} Primary 60G15; secondary 60G70

\section{Introduction}
Let $X(t),t\geq 0 $ be a centered stationary Gaussian process with
continuous trajectories, unit variance and correlation function $r$ satisfying for some $\alpha \in (0,2]$
$$1- r(t) \sim \abs{t}^{\alpha}, \quad t \rightarrow 0, \quad \text{ and } r(t)<1 , \quad \forall t> 0.$$
We write $\sim$ for asymptotic equivalence when the argument tends to 0 or infinity. \\
The seminal paper \cite{PicandsA} established for any $T$ positive  and $q(u)=u^{-2/\alpha}$
\BQN\label{pic}
\pk{\sup_{t\in[0,T]}X(t)>u}\sim T \mathcal{H}_{\alpha}     \frac{\pk{X(0)>u}}{q(u)} 
\EQN
as $u\to \IF$, where $ \mathcal{H}_{\alpha}$ is the {\it Pickands constant} defined by
\BQNY\label{pick}
\Ha=\lim_{T\rightarrow\infty}\frac{1}{T}\Ha[0, T] \in (0,\IF),\ \ \text{with}\ \ \Ha[0, T]=\E{
\sup_{t\in[0,T]}e^{ \sqrt{2}B_\alpha(t)-t^{\alpha}} },
 \EQNY
with $B_\alpha$ a standard fractional Brownian motion with Hurst index $\alpha/2$; see the recent contributions \cite{DiekerY, DM,MR3474465, SBK, KE17} for the main properties of  Pickands and related constants.

While the original proof of Pickands utilizes a discretisation approach,
in \cite{Pit72, Pit96} the asymptotics \eqref{pic} was derived
by establishing first  the exact asymptotics on the short interval $[0, q(u)T]$,
namely (see e.g., Lemma 6.1 in \cite{Pit96})
\BQN\label{eq:pick1}
\pk{\sup_{t\in[0, q(u) T]}X(t)>u}\sim \Ha[0, T] \pk{X(0)>u}, \quad u\to \IF
\EQN
and then using the {\it double-sum method}. A completely \KD{independent proof for the stationary case,
based on the notion of {\it sojourn time},
was derived by} Berman (see \cite{Berman82,Berman92}).

{In this contribution
we develop \kd{the {\it uniform double-sum method}. Originally, introduced by Piterbarg for non-stationary case, see e.g., \cite{Pit96},
the {\it double-sum method} is a powerful tool in derivation of the exact asymptotics}
of the tail distribution of supremum for non-stationary Gaussian processes (and fields).
With no \K{loss} of generality, for a given centered Gaussian process $Y(t), t\in [0,S]$
with continuous trajectories, the crucial steps of this method  are:\\
a) application of Slepian inequality that allows for uniform approximation as $u\to \IF$
(uniformity is with respect to $k \le N(u)$) of  summands of 
$\pk{\sup_{t\in[kTq(u) ,(k+1)Tq(u)]}Y(t)>u}$ by
$\pk{\sup_{t\in[0 ,Tq(u)]}X^\epsilon(t)>u_k}\kd{=:}p(u_k)$,
{\kd{for} appropriately chosen stationary process $X^\ve , \ve>0$;\\
b) uniform approximation for $k\le N(u)$ of $p(u_k)$ as $u\to \IF$;}\\
c) \kd{uniformly} tight upper bounds for the probability of double supremum
\begin{eqnarray}\label{double.1}
\pk{\sup_{t\in[kTq(u) ,(k+1)Tq(u)]}Y(t)>u,\sup_{t\in[lTq(u) ,(l+1)Tq(u)]}Y(t)>u}
\end{eqnarray}
for $k,l\in \mathcal{A}_u$, where the set $\mathcal{A}_u$ is suitably chosen.

The deep contribution \cite{dieker2005extremes} showed that
while dealing with supremum of Gaussian processes on the half-line
it is convenient to replace Slepian inequality by
a uniform version of the tail asymptotics of threshold-dependent Gaussian processes.
Omitting technical details,  \cite{dieker2005extremes}  
derives the exact asymptotics and a uniform upper bound of
$$\pk{ \sup_{t\in [0,T]} \xitu(t)> g_{u,\tau_u}}$$
 as $u \to \IF$, with respect to $\tau_u\in K_u$, for
 $\xitu$
\KD{being } centered Gaussian processes indexed by $u$ and $\tau_u$, see also Lemma 5.1 in \cite{KP2015}. This uniform counterpart of (\ref{eq:pick1}) is
crucial when the processes $X_{u,\tau_u}$
are parameterised by $u$ and $\tau_u$.

Recent contributions show \kd{strong need for analysis of distributional properties of more general
continuous functionals than supremum}, as  e.g.,
$\sup_{t\in [0, T]} \inf _{s\in [0,S]} X(s+f(u)t), S>0$, see \cite{MR3414985, MR3457055}
or $\inf_{s\in \mathcal{A}_u}\sup_{t\in \mathcal{B}_u} Y(s,t)$, see \cite{dekos14, KP2015}.\\
The lack of Slepian-type results for general continuous functionals $\Gamma$ can be overcome by the derivation of \rd{uniform} approximations with respect to $\tau_u$
of the tail distribution of $\Gamma (\xitu)$ as $u\to\infty$. Therefore, the principal goal of this contribution is to derive uniform approximations for
the tail of homogeneous \eM{continuous} functionals $\Gamma$ of general Gaussian random fields. Specifically,  we shall consider  $\Gamma$  defined on \rrd{$C(\bD)$}, the space of continuous functions on $E$ with $E \subset \R^d, d\ge 1$
a compact set containing the origin.   In Theorem \ref{lem51}
we derive the following uniform asymptotics
\BQN\label{eq1}
\limit{u} \sup_{\tau_u \in K_u} \ABs{
 \frac{ \pk{ \Gamma( \xitu)> \gtu} }{\Psi(g_{u,\tau_u})}  - C }=0,
\EQN
where $\xi_{u,\tau_u}(t), t\in E, \tau_u \in \KD{K_u}$ is a centered Gaussian random field, $C$ is a positive finite constant,  and  $\Psi$ denotes the survival function of an $N(0,1)$ random variable.  
This result 
 allows us to derive counterparts of (\ref{pic}) for a class of homogeneous
functionals of  centered Gaussian fields satisfying some weak asymptotic conditions. 
Additionally, in Section \ref{s.double} we derive a uniform upper bound for the double maxima
for general Gaussian fields parameterised  by
$u$ and $\tau_u$. That  extends and unifies the known upper bounds for (\ref{double.1}).

Brief organisation of the rest of the paper: \kd{main results of this contribution} and related discussions are presented in Section 2.
 We dedicate Section 3 to applications.
 Finally, we display the proofs of all the results in Section \ref{secprof}, \kd{postponing} some technical calculations to Appendix.

\section{Main Result}\label{s.main}
\kd{We begin this section with some motivations for the investigation
of distributional properties of functionals of threshold-dependent Gaussian random fields.
For this purpose we focus on supremum of non-centered Gaussian process.
Then we introduce the class of functionals that are of our interest and provide
the main result of this contribution; see Theorem \ref{lem51}.}

Numerous articles, e.g., \cite{HP99,debicki2002ruin, HP2004, dieker2005extremes},
developed techniques for the approximation, as $u\rw\IF$,
of the so-called ruin probability
\BQN
p(u)=\pk{\sup_{t \in \mathcal{T} } (X(t)-ct) >u},\label{ruin}
\EQN
where $X$ is a centered continuous Gaussian \kd{process}, $c>0$ is some constant and $\mathcal{T}= [0,\IF)$ or $\mathcal{T}=[0,T],T>0$.
\kd{Originally the {\it double-sum method} was designed to handle supremum
of centered Gaussian processes. For our case, this method still works under the following
\kk{modifications.} First, we rewrite the original problem in the language of a centered, threshold-dependent family of Gaussian processes
$Z_u(t)=\frac{X(t)}{u+ct}$, $u>0$ as follows
\BQN
p(u)  &=& \pk{\sup_{ t\in \mathcal{T} }Z_u(t)>1}.
\EQN
Then, one checks that, for suitably chosen $w(u)$ \kk{and $N(u)$,}
\BQN
p(u)
&\sim&
\pk{\text{There exists } {|k|\leq N(u)}: \sup_{t\in [0,w(u)S]}Z_u(t+kSw(u))>1}\nonumber\\
&\sim&
\sum_{|k|\leq N(u)}\pk{\sup_{t\in [0,S]} {Y}_{u,k}(t)>v_k(u)}=:\sum_{|k|\leq N(u)} p_k(u)
\label{motivation}
\EQN
as $u\rw\IF$ and $S\rw\IF$ respectively, where
$${Y}_{u,k}(t)=Z_u(w(u)t+w(u)kS)v_k(u), \quad v_k(u)=\inf_{t\in [0,S]}\frac{1}{\sqrt{Var(Z_u(w(u)t+w(u)kS))}}.$$
Finally, since usually
$ \limit{u} N(u)= \IF$, then
in order to determine the asymptotics
of $p(u)$ it is necessary to derive the asymptotics of
$p_k(u)$, as $u\to \IF$, uniformly for $\abs{k}\le N(u)$.
}
}

In this section, we consider a more general situation
\kd{focusing  on the validity of \eqref{eq1} for} centered Gaussian random fields.

\kd{
Next, let \ehd{$\bD \subset \R^d$ be a compact set including the origin and write $C(\bD)$ for the set of real-valued continuous functions defined on $\bD$}. Let $\Gamma: C(\bD)\to \R$   be a real-valued \eM{continuous} functional \ehd{satisfying}  \\
{\bf F1:}  there exists $c>0$ such that  $\Gamma(f)\le c \sup_{{t}\in\bD} {f({t})}$ for any \rrd{ $f\in C(\bD)$};\\
{\bf F2:}  
\K{$\Gamma(af+b)=a\Gamma(f)+b$ for any \rrd{$f\in C(\bD)$} and $a>0,b\inr$.}\\
{Note that {\bf F1-F2} cover the following important examples:
$$\Gamma=\sup, \quad \inf,  \quad a\sup+(1-a)\inf, \quad {a\in \mathbb{R}}.$$}
}
\kd{We} shall consider a family of  centered Gaussian random \k{fields}
$\xitu$ given by
$$\xitu(t)= \frac{Z_{u,\tau_u}(t)}{1+ h_{u,\tau_u}(t)}, \quad t\in E, \tau_u \in K_u,$$
with  $Z_{u,\tau_u}$ a centered Gaussian random field with unit variance and continuous trajectories,
and $h_{u,\tau_u}\in C_0(\bD)$, \ehd{where} $\coe$  {is the Banach space} of all continuous functions $f$ on $E$ such that $f(0)=0$ equipped with the sup-norm. In order to avoid trivialities, the thresholds $\gtu$ will be chosen such that
$$\limit{u} \pk{ \Gamma( \xitu)> \gtu}=0.$$
\rd{In order to derive the asymptotics of $\pk{\Gamma(\xitu)>\gtu }$ as $u\rw\IF$  we shall first condition on $\xitu({0}) = g_{u,\tau_u}-\frac{w}{g_{u,\tau_u}}$, yielding that
\BQNY
\pk{\Gamma(\xitu)>\gtu }
&=&
 \frac{ e^{- \gtu^2/2}}{\sqrt{2\pi}\gtu }
\int_{\R} e^{w-\frac{w^2}{2\gtu^2}}\pk{\Gamma(\chitu)> w}\, dw,
\EQNY
where
$$\chitu(t)=\gtu(\xitu({t})-\gtu)+w\Big|\Bigl( \xitu({0}) = g_{u,\tau_u}-\frac{w}{g_{u,\tau_u}}\Bigr).$$
Note that
\BQNY
\chitu(t) \KD{\stackrel{d}{=}} \frac{\gtu}{1+h_{u,\rd{\tau_u}}(t)}\Bigl(Z_{u,\tau_u}({t})-r_{u,\tau_u}({t},{0})Z_{u,\tau_u}({0})
\Bigr)+\EE{\chitu({t})}, \quad {t}\in\bD,
\EQNY
where $\KD{\stackrel{d}{=}}$ means equality of \kk{distributions}.}\\
Next,  we shall impose the following assumptions (\rdd{see also \cite{KP2015}[Lemma 5.1] and  \cite{dieker2005extremes}[Lemma 2]}) to ensure the weak convergence of  $\{\chitu(t), t\in E\}$, as $u\rw\IF$.   \\
{\bf C0}: The positive constants $\gtu$ are such that $ \limit{u} \inf_{  \gE{\tau_u}\in K_u } \gtu=\IF$. \\
{\bf C1}: \gE{There exists $h\in C_0(E)$ such that}
\BQN\label{PI1}
 \lim_{u\to \IF} \sup_{\tau_u\in K_u, t\in E} \abs{ \gtu^2 h_{u,\tau_u}(t)- h(t)}=0.
 \EQN
{\bf C2}: \KD{There exists
$\theta_{u,\tau_u}(s,t)$ such that
\BQN\label{PI2}
\lim_{u\rw\IF}\sup_{{\tau_u} \in K_u}\sup_{s\neq t\in E}\left| \gtu ^2\frac{Var\left(Z_{u,\tau_u}(t)-Z_{u,\tau_u}(s)\right)}
{2\theta_{u,{\tau_u}}(s,t)}-1\right| \eHHH{=} 0
\EQN
and for some centered Gaussian random field $\eta(t),t\inr^d$ with continuous trajectories
and $\eta(0)=0$
\BQN\label{A1}
\limit{u}\sup_{\tau_u \in K_u}|\theta_{u, \tau_u}(s, t)-Var(\eta(t)-\eta(s))|= 0, \quad \forall s, t\in E.
\EQN
}
{\bf C3}: There exists $a >0$  \KD{such that}
\BQN\label{A31}
\limsup_{u\rw\IF}\sup_{{\tau_u} \in K_u}\sup_{s\neq t, s,t\in E}\frac{\theta_{u, \tau_u }(s,t)}{\sum_{i=1}^{d}|s_i-t_i|^{a}}<\IF
\EQN
and
\BQN\label{eq2}
\lim_{\epsilon\downarrow  0}\limsup_{u\rw\IF}
\sup_{{\tau_u} \in K_u}\sup_{\norm{t-s}<\epsilon,s,t\in E}g_{u,{\tau_u}}^2
\EE{ \left[ Z_{u,{\tau_u}}(t)-Z_{u,{\tau_u}}(s)\right]Z_{u,{\tau_u}}(0)}=0.
\EQN
{If X is a centered Gaussian process with stationary increments sastifying
{\bf AI-AII} in \cite{KP2015}, then ${Y}_{u,k}(t), t\in [0,S], |k|\leq N(u)$ in
(\ref{motivation}) satisfies {\bf C0-C3};
\kk{see also \cite{dieker2005extremes}.} \\
   The  intuitive explanation behind these assumptions is as follows:
\ehd{   {\bf C1} and (\ref{eq2}) in {\bf C3} are used to guarantee the
uniform convergence of the function $\EE{\chitu({t})}$ for $t\in E$ as $u\rw\IF$.
Utilising further {\bf C2}, the convergence of finite-dimensional distributions
\ehd{(fidi's)} of $\chitu({t}), t\in E$  to those of $\eta(t),t\in E$ can be shown.
Moreover, the tightness follows by (\ref{A31}) in {\bf C3}.}  \COM{$\{g_{u,\tau_u}\Big(Z_{u,\tau_u}({t})-r_{u,\tau_u}({t},{0})Z_{u,\tau_u}({0})\Big), t\in E\}$ to $\{\eta(t), t\in E\}$ and (\ref{A31}) in {\bf C3} shows the tightness of the process $\left\{g_{u,\tau_u}\Bigl(Z_{u,\tau_u}({t})-r_{u,\tau_u}({t},{0})Z_{u,\tau_u}({0})\Bigr), t\in E\right\}, u>0$. Therefore the weak convergence is satisfied.
}

\def\Hh{\mathcal{H}_{\eta, h}^{\Gamma}}

Given $h \in  \coe $ and  \ehd{the functional} $\Gamma$ satisfying {\bf F1-F2}, \kd{for $\eta$ introduced in {\bf C2},} we define a new constant
\BQN\label{eq:GPitCons}
{\Hh(E):=\EE{e^{\Gamma(\eta^h )}}, \quad {\eta^h(t):= \sqrt{2}\eta(t)-Var(\eta(t))-h(t)}, }
\EQN
which by {\bf F1} is finite.  For notational simplicity we set below
$$\mathcal{H}_\eta(E)=\mathcal{H}_{\eta, 0}^{\sup}(E).$$
We present next \kd{the main result of this section}.
\KD{Recall that
$\Psi$ stands for the survival function of an $N(0,1)$ random variable.}

\BT \label{lem51} Under assumptions {\bf C0-C3} and {\bf F1-F2}, if further $\pk{ \Gamma( \xitu)> \gtu} >0$ for all $\tau_u \in K_u$ and all $u$ large, then
\BQN\label{eq1_m}
\limit{u} \sup_{\tau_u \in K_u} \ABs{
 \frac{\pk{ \Gamma( \xitu)> \gtu}}{\Psi(g_{u,\tau_u})}  - \Hh(E) }=0.
\EQN
\ET

\BRM\label{remark}
i) \K{Under} the assumptions of \netheo{lem51} we have
\BQN
 \limsup_{u\to \IF} \sup_{\tau_u \in K_u}  \frac{\pk{ \Gamma( \xitu)> \gtu}}{\Psi(\gtu)}   < \IF ,
 \EQN
{which coincides with the results of Lemma 5.1 in \cite{KP2015} and extends Lemma 2 in \cite{dieker2005extremes}.}\\
{ii) Condition} {\bf C2} and (\ref{eq2}) in {\bf C3} are equivalent to {\bf C2} and
\BQN\label{eq3}
\lim_{u\rw\IF} \sup_{t\in E, {\tau_u} \in K_u} \ABs{ g_{u,\tau_u}^2Var(Z_{u,\tau_u}({t})-Z_{u,\tau_u}({0}))- 2 Var(\eta({t}))}=0.
\EQN
 iii) Condition {\bf C2} can be formulated also for the degenerated case $\eta(t)=0,t\inr^d$ almost surely. The claim of \netheo{lem51} holds also for such $\eta$.
\ERM
\def\qq{q}
\def\qqq{b}

Next we give a simplified version of Theorem \ref{lem51}.  Instead of {\bf C2-C3}, we  assume that \BQN\label{A2}
\lim_{u\rw\IF}\sup_{{\tau_u} \in K_u}\sup_{s\neq t, s,t\in E}\left| \gtu ^2\frac{Var\left(Z_{u,\tau_u}(t)-Z_{u,\tau_u}(s)\right)}
{2\sum_{i=1}^d \frac{c_i\sigma_i^2(\qq_i(u)|s_i-t_i|)}{\sigma_i^2(\qq_i(u))}}-1\right| \eHHH{=} 0,
\EQN
where
$\qq_i(u), i=1,\dots, d$ are some functions of $u$ with $\qq_i(u)>0$ for $u$ large enough and $\lim_{u\rw\IF} \qq_i(u)=\varphi_i\in [0,\IF]$ with $$\varphi_i=\left\{\begin{array}{cc}
0, & 1\leq i\leq d_1\\
(0,\IF), & d_1+1\leq i\leq d_2,\\
\IF, & d_2+1\leq i\leq d
\end{array}\right. $$
and  $\ehd{c_i\ge 0}, 1\leq i\leq d$.
Moreover, \rd{ $\sigma_i, 1\leq i\leq d$} are   regularly varying at $0$ with indices
$\alpha_{i,0}/2\in (0,1]$ respectively and $\sigma_i(0)=0$, $\sigma_i(t)>0, t>0$, $1\leq i\leq d$;
 $\sigma_i, d_2+1\leq i\leq d$ are  bounded on any compact interval and regularly varying at
$\IF$ with indices $\alpha_{i,\IF}/2\in (0,1]$, respectively;  $\sigma_i^2(t), d_1+1\leq i\leq d_2$ are \rrd{continuous} and  non-negative definite, implying that there exist centered Gaussian processes  $\eta_i, \rd{d_1+1\leq i\le d_2}$  with continuous sample path and  \rd{stationary increments such that  }
 $Var(\eta_i(t)):=\sigma_i^2(t), d_1+1\leq i\leq d_2$.
 \rdd{
We refer to, e.g., \cite{debicki2002ruin,dieker2005extremes,HP99,HP2004}, where particular examples of Gaussian processes that
satisfy the above regularity assumptions are investigated; see also \cite{K2010}  for characterisation of such processes in terms of max-stable stationary processes.}\\

\BS\label{ths} Suppose that {\bf C0-C1} and {\bf F1-F2} hold.
If  \kk{(\ref{A2})
holds
with $\sum_{i=1}^d c_i>0$ and
$\pk{ \Gamma( \xitu)> \gtu} >0$ for all $\tau_u \in K_u$ and all $u$ large, then}  (\ref{eq1_m}) holds with
\BQN\label{A3}
\eta(t)=\sum_{i=1}^{d_1} \sqrt{c_i}B_{\alpha_{i,0}}(t_i)+\sum_{i=d_1+1}^{d_2} \sqrt{c_i}\frac{\eta_i(\varphi_it_i)}{\sigma_i(\varphi_i)}+\sum_{i=d_2+1}^{d} \sqrt{c_i}B_{\alpha_{i,\IF}}(t_i),
\EQN
\rrd{where $B_{\alpha_{i,0}}, 1\leq i\leq d_1$, $\eta_i, d_1+1\leq d_2$ and $B_{\alpha_{i,\IF}}, d_2+1\leq i\leq d$ are mutually independent.}
\ES

\BRM i)
\kd{Condition (\ref{A2}) is satisfied by a large class of important processes that are investigated in the literature,} see e.g.
\cite{HP99,debicki2002ruin, dieker2005extremes,KP2015,KEP2016}.\\
ii) Under the assumptions of \netheo{lem51}
\BQN\label{grosh}
\limit{u} \sup_{\tau_u \in K_u} \ABs{ \frac{ \pk{ \Gamma_i( \xitu) > u, i=1 \ldot d }}{\Psi(\gtu)}  -
	\mathcal{H}^{ \Gamma_1\ldot \Gamma_d}_{\eta, h}}=0,
\EQN
with $\Gamma_i, i\le d $ continuous functionals satisfying {\bf F1-F2} and
$$ \mathcal{H}^{ \Gamma_1\ldot \Gamma_d}_{\eta, h}=\int_{\R} e^w \pk{ \Gamma_i(\ehd{\eta^h})> w, i = 1 \ldot d}
\, dw\in (0,\IF).$$
Moreover, \eqref{grosh} holds also in the case that $\eta$ is degenerated, i.e.,  $\eta(t)=0,t\inr^d$ almost surely.
\ERM

\kd{Finally, we present below a version
of \netheo{lem51} under slightly different and more
explicit assumptions. We keep the same notation as in Theorem \ref{lem51}} \rrd{and moreover let
$\sigma_{u,\tau_u}^2(t):= Var( \xitu(t))$}.  \\
{\bf D1:} Condition {\bf C0} holds for  $g_{u,{\tau_u}}$ and 
$\sigxiu({0})=1$ for all  $\tau_u \in K_u$ and all $u>0$, and there exists some  $h \in \coe$  such that
$$\lim_{u\rw\IF}\sup_{{t}\in\bD, \gE{\tau_u} \in K_u}\abs{\gtu^2(1-\sigxiu({t}))-h({t})}=0.$$
{\bf D2:} There exists a centered Gaussian random field $\eta({t}), {t}\in\R^d$ with continuous sample paths, $\eta({0})=0$ such that for any  ${s},{t}\in  \bD$ and ${\tau_u}\in K_u$
\BQN \label{c2a}
\lim_{u\rw\IF} \sup_{{\tau_u} \in K_u} \ABs{ \gtu^2Var(\xitu({t})-\xitu({s}))- 2 Var(\eta({t})-\eta({s}))}=0,
\EQN
and
\BQN \label{c2b}
\lim_{u\rw\IF} \sup_{t\in E, \gE{\tau_u} \in K_u} \ABs{ \gtu^2Var(\xitu({t})-\xitu({0}))- 2 Var(\eta({t}))}=0.
\EQN
{\bf D3:} There exist positive constants $G, \nu, u_0$ such that for any $u>u_0$  
\BQNY
\sup_{\tau_u \in K_u}  \gtu^2Var(\xitu({t})-\xitu({s})) \le G \norm{{t}-{s}}^{\nu}
\EQNY
holds for all ${s},t\in  \bD$.

\BT\label{cor} If {\bf D1-D3} and {\bf F1-F2} are satisfied, then (\ref{eq1_m}) holds.
\ET

\section{Applications}
\subsection{Upper Bounds for Double Supremum}\label{s.double}
\KD{Uniform bounds for the tail distribution of bivariate maxima of Gaussian processes
play a key role in the double-sum technique of V.I. Piterbarg; see, e.g., \cite{Pit96,Pit20}.}
More precisely, of interest is to find an optimal upper bound for
$$
D(\lambda_1, \lambda_2, \mathcal{E}_1, \mathcal{E}_2, u):= \pk{ \sup_{t\in \lambda_1+\mathcal{E}_1} X_u(t)> m_{\lambda_1}(u), \sup_{t\in\lambda_2 +\mathcal{E}_2} X_u(t)> m_{\lambda_2}(u)},
$$
which is valid for all large $u$ with  $\lambda_i$'s and $\mathcal{E}_i$'s  controlled by $E_u$ by requiring that $\lambda_i + \mathcal{E}_i \subset E_u$, \ehd{with $E_u$ a compact subset of $\R^d$}.
Further, the thresholds $m_{\lambda_1}(u), m_{\lambda_2}(u)$  are assumed to satisfy
  \BQN\label{horizon}
  \lim_{u\rw\IF}m(u)=\IF,\ \  \lim_{u\rw\IF}\sup_{\lambda_i + \mathcal{E}_i \subset E_u}\left|\frac{m_{\lambda_i}(u)}{m(u)}-1\right|=0, \quad  i=1,2
  \EQN
for some positive function $m$.

Set below  $F(A,B)= \inf_{s\in A, t\in B}\norm{s-t}$ with  $A, B$ \eM{two non-empty subsets of} $\mathbb{R}^d$ and $\norm{\cdot}$ the Euclidean norm.
Let $\mathbb{K}=\{{(\lambda_1, \lambda_2):}\lambda_i+\mathcal{E}_i\subset E_u, i=1,2\}$.

\BT\label{pro}  Let $X_u(t), t\in E_u\subset \R^d$ be a family of  centered Gaussian random \ehd{fields} with continuous trajectories, variance $1$ and correlation function $r_u$.
 Suppose that there exist positive constants $S_1,\mathcal{C}_1, \mathcal{C}_2, \beta$ and $\alpha \in (0,2]$  such that for  $u$ sufficiently large
\BQN\label{fu2}
 &&m^2(u)(1-r_u(s, t))\geq \mathcal{C}_{1}\norm{s-t}^{\beta},  \norm{s-t}\geq S_1,
 \quad s,t\in E_u,
 \EQN
 and
 \BQN\label{COR}
 m^2(u)(1-r_u(s, t))\leq \mathcal{C}_2 \norm{s-t}^\alpha, \ \  s,t\in E_u, s-t\in [-1, 1]^d.
 \EQN
 Moreover, there exists $\delta>0$ such that for $u$ large enough
 \BQN\label{delta}
 r_u(s,t)>\delta-1, \quad  s,t\in E_u.\EQN
 If further  (\ref{horizon}) holds, then
there exists $ \mathcal{C}>0$ such that  for \ehd{all} $u$ large enough
 \BQN\label{unibound}
  \sup_{(\lambda_1,\lambda_2)\in \mathbb{K}, \mathcal{E}_i\subset [0,S_2]^d, \mathcal{E}_i\neq \emptyset, i=1,2 }     \frac{ e^{\frac{\mathcal{C}_1 F^\beta(\lambda_1+\mathcal{E}_1, \lambda_2+\mathcal{E}_2)}{8}}D(\lambda_1, \lambda_2, \mathcal{E}_1, \mathcal{E}_2, u)}{ \HEH{S_2}^{2d} \Psi(m_{\lambda_1, \lambda_2}(u)) }
& \le & \mathcal{C} ,
 \EQN
with $S_2>1$, $m_{\lambda_1, \lambda_2}(u)=\min(m_{\lambda_1}(u), m_{\lambda_2}(u))$ and $\mathcal{C}$ a positive constant independent of $S_2,u$.
 \ET
 \cl{Next assume that  $\kappa_i(t)>0, t>0, 1\leq i\leq 2d$ are some non-negative  locally bounded functions and \ehd{define}
 $$g_u(s,t)=\sum_{i=1}^{d}\frac{\kappa_i(\qq_i(u)|s_i-t_i|)}{\kappa_i(\qq_i(u))} \quad \text{and}
 \quad \widetilde{g}_u(s,t)=\sum_{i=1}^{d}\frac{\kappa_{i+d}(\qq_{i+d}(u)|s_i-t_i|)}{\kappa_{i+d}(\qq_{i+d}(u))}.$$
Further, let  $\qq_i(u)>0, u>0$ be such that  $$\lim_{u\rw\IF}\qq_i(u)=\varphi_i\in [0,\IF],
 \quad \ehd{1 \le i \le 2d.}$$
\BK\label{CC0} Let $X_u(t), t\in E_u$ be  centered Gaussian random fields with continuous trajectories, variance $1$ and correlation function $r_u$
 satisfying (\ref{delta}). Assume further that (\ref{horizon}) holds. If further
  for $u$ sufficiently large
  \BQN\label{corr1}
  \mathcal{C}_3 \k{g_u}(s,t) \leq  m^2(u)(1-r_u(s,t))\leq \mathcal{C}_4 \k{\widetilde{g}_u}(s,t), \quad s,t\in E_u,
  \EQN
  with $\mathcal{C}_3, \mathcal{C}_4>0$ and $\kappa_i, 1\leq i\leq 2d$, \HEH{being
 regularly varying} both at $0$ and \K{at} $\IF$ with indices $\alpha_{i,0}>0$ and $ \alpha_{i,\IF}>0$, respectively,
    then there exists $ \mathcal{C}>0$ such that for $u$ large enough (\ref{unibound}) holds with $\beta=\frac{1}{2}\min_{i=1,\dots, 2d}\min(\alpha_{i,0}, \alpha_{i,\IF}, 2)$ and $\mathcal{C}_1$ a fixed positive constant.
 \EK
\BK\label{CC1} Let $X_u(t), t\in E_u\ehd{\subset \R^d}$ be  centered Gaussian random fields with continuous trajectories, variance $1$ and correlation function $r_u$
 satisfying (\ref{delta}) and (\ref{corr1}) with $\varphi_i=0, 1\leq i\leq 2d$ and $\kappa_i, 1\leq i\leq 2d$ \HEH{being
 regularly varying}  at $0$ with indices $\alpha_{i,0}>0$. If  further (\ref{horizon}) and
 \BQN\label{free}
 {\limsup_{u\rw\IF}\sup_{s, t\in E_u}\ehd{\max}_{i=1,\dots, 2d}\qq_i(u)|s_i-t_i|<\IF}
 \EQN
hold, \K{then} there exist positive constants $ \mathcal{C}, \mathcal{C}_1$ such that for $u$ large enough (\ref{unibound}) holds with $\beta=\frac{1}{2}\min(2, \min_{i=1,\dots, 2d} \alpha_{i,0})$.
\EK
}

\BRM i) Under the assumptions of Theorem \ref{pro}, {using the idea of \cite{MR2145669,DebickiMandj10},
since for $\gamma\in (0,1)$}
\BQNY
D(\lambda_1, \lambda_2, \mathcal{E}_1, \mathcal{E}_2, u)
&\leq& \pk{ \sup_{s\in \lambda_1+\mathcal{E}_1, t\in\lambda_2 +\mathcal{E}_2}( \gamma X_u(s)+(1-\gamma)X_u(t))>  m_{\lambda_1, \lambda_2, \gamma}(u)},
\EQNY
\k{with  $ m_{\lambda_1, \lambda_2, \gamma}(u)
= \gamma m_{\lambda_1}(u)+ (1-\gamma)m_{\lambda_2}(u)$,
}
\k{then in some cases} \eqref{unibound} can be improved by putting
\k{$4\gamma (1- \gamma) \mathcal{C}_1$}
instead of
$\mathcal{C}_1$ and
$m_{\lambda_1, \lambda_2, \gamma}(u)$
instead of $m_{\lambda_1, \lambda_2}(u)$, respectively. \\
ii) A particular example  is $\kappa_i(x)= x^{\alpha_i}, \alpha_i \in (0,2]$. For such a case, the result of Corollary \ref{CC1} yields the claim of Lemma 9.14 in
\cite{Pit20}, see also Lemma 6.3 in \cite{Pit96}.\\
\ERM

\subsection{\HEH{Tail Approximation of  $\Gamma_{E_u}(X_u)$}}\label{s.gen.funct}

In many applications the tail asymptotics of general functionals of Gaussian random fields
$X_u$  indexed by thresholds $u>0$ is of interest. In this section we present an
application
\K{of Theorem \ref{lem51}} concerned with the tail asymptotics of $\Gamma_{E_u}(X_u)$, where
$$E_u:=\left(\prod_{i=1}^d[a_i(u), b_i(u)]\right)\times E$$
is also parametrised by $u$, with $E$ a compact subset of $ \mathbb{R}^n,n\inn$. Without loss of generality, we assume $0\in E$.
The functional $\Gamma_{E_u}$ is defined as follows:\\
  Let $\Gamma^*: C(E)\to \R$   be a real-valued \eM{continuous} functional satisfying {\bf F1-F2} with $c=1$ in {\bf F1}. For any  compact \ehd{set} $A\subset \mathbb{R}^d$ define
$$\Gamma_{A\times E}(f)=\sup_{s\in A}\Gamma^*(f(s,t)), \quad f\in C(A\times E).$$
It follows that $\Gamma_{A\times E}$ is \HEH{a continuous functional} and satisfies {\bf F1--F2} with $c=1$ in {\bf F1}. Examples of $\Gamma^*$ are
$$\Gamma^*=\sup, \quad \inf,  \quad a\sup+(1-a)\inf, \quad {a\leq 1}.
$$
We shall consider  $X_u(s,t), (s,t)\in E_u$, a family of centered continuous
Gaussian random  fields with variance function $\sigma_u(s,t)$ and correlation function
$r_u(s,t,s',t')$ satisfying {as $u\rw\IF$}
\rdd{\BQN\label{var2}
\sigma_u(0,0)=1, \ \ 1-\sigma_u(s,0)\sim \sum_{i=1}^d\frac{|s_i|^{\beta_i}}{g_i(u)}, \quad  s\in \prod_{i=1}^d[a_i(u), b_i(u)]
\EQN
and
\BQN\label{new1}
\lim_{u\rw\IF}\sup_{s\in \prod_{i=1}^d[a_i(u), b_i(u)], t\neq 0, t\in E}\left|\frac{1-\frac{\sigma_u(s,t)}{\sigma_u(s,0)}}{\sum_{i=d+1}^{d+n}\frac{|t_i|^{\beta_i}}{g_i(u)}}-1\right|=0,
\EQN
}
where $\beta_i>0$ and  $g_i(u)$ is a function of $u$ satisfying $\lim_{u\rw\IF}g_i(u)=\IF$ for $1\leq i\leq d+n.$
Moreover, \kd{there exists $m(u)$ such that $\lim_{u\rw\IF}m(u)=\IF $ and}
\BQN\label{cor2}
\lim_{u\rw\IF}\sup_{(s,t), (s',t')\in E_u, (s,t)\neq (s',t')}\left| \frac{\ehd{m^2(u)}( 1-r_u(s,t,s',t'))}{\sum_{i=1}^d{\frac{c_i\sigma_i^2(\qq_i(u)|s_i-s_i'|)}{\sigma^2_i(\qq_i(u))}}+\sum_{i=d+1}^{d+n}
\frac{c_i\sigma_i^2(\qq_i(u)|t_i-t_i'|)}{\sigma_i^2(\qq_i(u))}}-1\right|=0,
\EQN
where  $ c_i>0$,  $\qq_i(u)>0$, $\lim_{u\rw\IF}\qq_i(u)=\varphi_i\in [0,\IF],  1\leq i\leq d+n$,
\k{and $\sigma_i$ are the variance functions of
$\eta_i$'s, centered continuous Gaussian processes with  stationary increments, \HEH{$\eta_i(0)=0$},
satisfying \ehd{further the following assumptions:}\\
{\bf A1:} $\sigma_i^2(t)$ is regularly varying at $\IF$ with index $2\alpha_{i,\IF}\in (0,2)$
and is continuously differentiable over $(0,\IF)$ with $\dot{\sigma_i^2}(t)$ being ultimately monotone at $\IF$.\\
{\bf A2:} $\sigma_i^2(t)$ is regularly varying at $0$ with index $2\alpha_{i,0}\in (0,2]$.}\\
Moreover, we shall assume \kd{that}
$$\lim_{u\rw\IF}\frac{|a_i(u)|^{\beta_i}}{g_i(u)}
=\lim_{u\rw\IF}\frac{|b_i(u)|^{\beta_i}}{g_i(u)}=0, \quad 1\leq i\le d\kd{+n}. $$
\cl{Let
\BQN\label{vf} V_{\varphi_i}(t_i)=\left\{\begin{array}{cc}
\sqrt{c_i}B_{\alpha_{i,0}}(t_i), & \varphi_i=0\\
\frac{\sqrt{c_i}}{\sigma_i(\varphi_i)}\eta_i(\varphi_it_i), & \varphi_i\in (0,\IF), \\
\sqrt{c_i}B_{\alpha_{i,\IF}}(t_i), & \varphi_i=\IF
\end{array}\right. \quad 1\leq i\leq d+n.
\EQN
}
{In the sequel, we shall denote
	$$\mathcal{P}_\eta^h(E)=\mathcal{H}_{\eta, h}^{\sup}(E),  \quad \mathcal{H}_\eta(E)=\mathcal{H}_{\eta, 0}^{\sup}(E)$$
and set
	$$\mathcal{P}_\eta^h=\lim_{S\rw\IF}\mathcal{P}_\eta^h([0,S]), \quad \widehat{\mathcal{P}}_\eta^h=\lim_{S\rw\IF}\mathcal{P}_\eta^h([-S,S]), \quad \mathcal{H}_\eta=\lim_{S\rw\IF} S^{-1} \mathcal{H}_\eta([0,S])$$
	if the limits exist. We refer to \cite{Pit96, MR1993262,  KEP2016}
for the properties of {\it Piterbarg constants} $\mathcal{P}_\eta^h$}
and {\it Pickands constants} $\mathcal{H}_\eta$. {Next, suppose that}
$$\lim_{u\rw\IF}\frac{m^2(u)}{g_i(u)}=\gamma_i\in [0,\IF]$$
and for all $u$ large
$\pk{\Gamma_{E_u}(X_u)>m(u)}>0$. \\
\BT\label{thm}
Let $X_u(s,t), (s,t)\in E_u\ehd{\subset \R^{d+n}}$ be a family of centered  Gaussian random fields with continuous
trajectories satisfying \kd{(\ref{var2})-(\ref{cor2})} and
 $$
 \HEH{\gamma_i = \left\{
  \begin{array}{cc}
0 ,&  \text{if} \ 1\leq i\leq d_1,\\
 \IF ,&  \text{if} \ d_2+1\leq i\leq d,\\
 \end{array} \right.
\quad \gamma_i \in (0,\IF),\ d_1+1\leq i\leq d_2, \quad \gamma_i\in [0,\IF),  \ d+1\leq i\leq d+n.
 }$$
If further for $1\leq i\leq d_1$ $$\lim_{u\rw\IF}\frac{(m(u))^{2/\beta_i}a_i(u)}{(g_i(u))^{1/\beta_i}}=y_{i,1}, \quad \lim_{u\rw\IF}\frac{(m(u))^{2/\beta_i}b_i(u)}{(g_i(u))^{1/\beta_i}}=y_{i,2},
\quad  \lim_{u\rw\IF}\frac{(m(u))^{2/\beta_i}(a_i^2(u)+b_i^2(u))}{(g_i(u))^{2/\beta_i}}=0,$$
 with $-\IF\leq y_{i,1}<y_{i,2}\leq \IF$, for $d_1+1\leq i\leq d_2,$  $ a_i(u)\leq 0\leq b_i(u), \lim_{u\rw\IF}a_i(u)=a_i\in[-\IF,0], \lim_{u\rw\IF}b_i(u)=b_i\in [0,\IF] $ and $a_i(u)\leq 0\leq b_i(u)$ for $d_2+1\leq i\leq d$,  then
\BQN\label{mainm}
\lefteqn{\pk{\Gamma_{E_u}(X_u)>m(u)}}\nonumber\\
&  \sim& \prod_{i=1}^{d_1}\mathcal{H}_{V_{\varphi_i}}\prod_{i=d_1+1}^{d_2}
\mathcal{P}_{V_{\varphi_i}}^{h_{i}}[a_i,b_i]\rd{\mathcal{H}_{\widetilde{V}_{\varphi }, \widetilde{h}}^{\Gamma^*}(\HEH{E})}\prod_{i=1}^{d_1}\int_{y_{i,1}}^{y_{i,2}}e^{-|s|^{\beta_i}}ds \prod_{i=1}^{d_1}\left(\frac{g_i(u)}{m^2(u)}\right)^{1/\beta_i}\Psi(m(u)),
\EQN
where
\BQN\label{vf1}
\widetilde{V}_{\varphi}(t)=\sum_{i=1}^{n}V_{\varphi_{d+i}}(t_i),  \quad \widetilde{h}(t)=\rd{\sum_{i=1}^{n}\gamma_{d+i}|t_i|^{\beta_{d+i}}},
 \quad h_{i}(s_i)=\gamma_i|s_i|^{\beta_i}, \quad d_1+1\leq i\leq  d_2.
 \EQN
\ET
\BRM
\K{Theorem \ref{thm} 
extends and unifies both the previous findings of \cite{HP99, debicki2002ruin, HP2004, dieker2005extremes}
\ehd{and  in particular Theorem 8.2 in}  \cite{Pit96}.}
\ERM

\subsection{Generalized Piterbarg Constants}
Let $X(t),t \geq 0$ be a centered Gaussian process with stationary increments and continuous trajectories. Suppose that the variance function $\sigma^2(t)=Var(X(t))$ is strictly positive for all $t>0$ and $\sigma(0)=0$.
Define next 
\bqny{
	\mathcal{P}_X^{b}(\ehd{[0,S], [0,T]})=\E{\eM{\sup_{t\in[0,T]}\inf_{ s\in[0,S]}}
		e^{ \sqrt{2}X(t-s)-(1+ b)\sigma^2(\abs{t-s})}},}
where $b, S,T$ are positive constants. In the special case, that $X=B_\alpha$
is a fractional Brownian motion (fBm) with Hurst index $\alpha/2 \in (0,1]$,
the generalized Piterbarg constant
$$\mathcal{P}_{B_\alpha}^{b}(S)= \limit{T} \mathcal{P}_{B_\alpha}^{b}(\ehd{[0,S], [0,T]})\in (0,\IF)$$
determines the asymptotics of Parisian ruin of the corresponding risk model, see \cite{MR3457055}.
\K{Note that the classical Piterbarg constant corresponds to the case $S=0$.}
Our next result shows that  $\mathcal{P}_X^{b}(S) \in (0,\IF)$ for a general Gaussian process  with stationary increments.

\BS \label{propP}
\k{If  $X(t),t \geq 0$	is a centred Gaussian	process with stationary increments and variance function
satisfying {\bf A1}  with regularly varying index $2\alpha_\IF\in (0,2]$ and {\bf A2} with regularly varying index $2\alpha_0\in (0,2)$,  then  for any $b,S$ positive we have}
\BQNY \label{boundPit}
\lim_{T\to \IF} \mathcal{P}_X^{b}([0,S], [0, T]) \in (0,\IF) .
\EQNY
\ES

\def\gt{u}
\def\xiu{\xi_u}

\section{Proofs}
\label{secprof}

\ehd{Hereafter,  by $\mathbb{Q}$, $\mathbb{Q}_i, i=1,2,\dots$ we denote positive constants which may differ from line to line.}\\

\prooftheo{lem51} Since we assume that $ \pk{\Gamma(\xitu)>\gtu }>0$
for all $u$ large and any $\tau_u \in K_u$, then  by \KD{conditioning} 
\BQNY
\pk{\Gamma(\xitu)>\gtu }&=& \int_{\R} \pk{\Gamma(\xitu)>\gtu  \lvert \xitu({0}) = x} \frac{e^{- x^2/2}}{\sqrt{2 \pi}}\, dx \\
&=&
 \frac{ e^{- \gtu^2/2}}{\sqrt{2\pi}\gtu }
\int_{\R} e^{w-\frac{w^2}{2\gtu^2}}\pk{\Gamma(\chitu)> w}\, dw\\
&=:& \frac{ e^{- \gtu^2/2}}{\sqrt{2\pi}\gtu }\mathcal{I}_{u,\tau_u},
\EQNY
with $\mathcal{I}_{u,\tau_u}>0$ \ehd{for all $u$ large} and
$$\chitu(t)=\zetatu (t)| (\zetatu({0})=0), \quad \zetatu(t)= \gtu(\xitu({t})-\gtu)+w.$$
Hence the proof follows by \ehd{showing that $\Hh(E)$ is finite and}
\BQN \label{esB}
\limit{u} \sup_{\tau_u \in K_u} \abs{\mathcal{I}_{u,\tau_u}- \Hh(E)} = 0.
\EQN
{\it \underline{Weak convergence of $\Gamma(\chitu)$}}.
We have that $\chitu(0)=0$ almost surely. Setting  $r_{u,\tau_u}({s},{t})=\rd{Cor}(Z_{u,\tau_u}({s}), Z_{u,\tau_u}({t}))$ we may write
\BQNY
\chitu(t) \KD{\stackrel{d}{=}} \frac{\gtu}{1+h_{u,\rd{\tau_u}}(t)}\Bigl(Z_{u,\tau_u}({t})-r_{u,\tau_u}({t},{0})Z_{u,\tau_u}({0})\Bigr)+\EE{\chitu({t})}, \quad {t}\in\bD,
\EQNY
where $\KD{\stackrel{d}{=}}$ means equality of the fidi's.   Since
$$
(1+h_{u,\tau_u}(t))\EE{\chitu({t})}=- \gtu^2(1-r_{u,\tau_u}({t},{0}))-g_{u,\tau_u}^2h_{u,\tau_u}(t)+w(1-r_{u,\tau_u}(t,0)+h_{u,\tau_u}(t))
$$
by {\bf C1, C3} for some arbitrary $M$ positive, uniformly with respect to $t\in\bD, \tau_u \in K_u, w\in [-M,M]$
\BQN\label{eqex}
(1+h_{u,\tau_u}(t))\EE{\chitu({t})}
& \to & -(\sigma_\eta^2({t})+h({t})), \  u\rw\IF
\EQN
and also for any $s,t \in \bD$ uniformly with respect to $ \tau_u \in K_u, w\in [-M,M]$
\BQN\label{eqco}
\lefteqn{Var\Bigl( (1+h_{u,\tau_u}(t))\chitu({t})-(1+h_{u,\tau_u}(s))\chitu({s})\Bigr)}\notag \\
&=&\gtu^2\LT[\EE{\Bigl( Z_{u,\tau_u}({t})-Z_{u,\tau_u}({s})\Bigr)^2}-\LT(\EE{Z_{u,\tau_u}({0})[Z_{u,\tau_u}({t})-Z_{u,\tau_u}({s})]}\RT)^2\RT]\notag \\
& \to & 2 Var(\eta({t})-\eta({s})), \quad  u\rw\IF.
\EQN
Consequently, by Lemma 4.1 in \cite{Yimin15} the \ehd{fidi's of}
 $ (1+h_{u,\tau_u}(t))\chitu(t), t\in E$ converge
 to those of \rd{$\eta^h(t), t\in E$ } \kd{as} $u\to \IF$
uniformly for $\tau_u \in K_u, w\in [- M,M]$ where
$M>0$ is fixed {(recall  $\eta^h(t)=\sqrt{2}\eta(t)-Var(\eta(t))-h(t)$)}. \\
Condition {\bf C3} together with the uniform convergence in
\eqref{eqex} guarantee that Proposition 9.7 in \cite{Pit20} can
be applied to yield the uniform tightness
of $(1+h_{u,\tau_u}(t))\chitu(t), t\in E$ and thus
$\{(1+h_{u,\tau_u}(t))\chitu(t), t\in E\}$ weakly \kd{converges} to
$\{\eta^h(t), t\in E\}$, as $u\to \IF$,  uniformly with respect to $\tau_u \in K_u$.
 Further, since
\ehd{$$ \limit{u} \sup_{ t \in E, \tau_u \in K_u} h_{u,\tau_u}(t) =0, $$}
then $\{\chitu(t), t\in E\}$ converges weakly  to $\{\eta^h(t), t\in E\}$ as $u\to \IF$,  uniformly with respect to $\tau_u \in K_u$.\\
Consequently, since \eM{we assume that} $\Gamma$ is a continuous functional, by the continuous mapping theorem $\Gamma(\chitu)$ converges in distribution to $\Gamma({\eta^h})$ as $u\to \IF$ uniformly with respect to \ehd{$\tau_u\in K_u$}. \\
{\it \underline{Convergence of (\ref{esB}})}.
\rdd{Denote by $\mathbb{A}=\{w: \pk{\Gamma(\eta^h )> \kd{w}} \text{ is discontinuous at  } w\}$, then $\mathbb{A}$ is an countable set with measure $0$.
Hence for any $w\in \mathbb{R}\setminus \mathbb{A}$
$$\limit{u}\sup_{\tau_u \in K_u\COM{w\in [-M,M]}}\ABs{\pk{\Gamma(\chitu)> w}- \pk{\Gamma(\eta^h )> w} }=0     $$
and by ${\bf C0}$
$$\limit{u}\sup_{\tau_u \in K_u, w\in [-M,M]}e^w \Bigl[1-   e^{-\frac{w^2}{2 \gtu^2}} \Bigr] \le
 \frac{e^M M^2}{2 \liminf_{u\to \IF }\inf_{\gE{\tau_u}\in K_u}  \gtu^2} \to 0, \quad u \to \IF      $$
implying
\BQNY
&&\limit{u}\sup_{\tau_u \in K_u}\ABs{\int_{-M}^M
\Biggl[ e^{w-\frac{w^2}{2 \gtu^2}}\pk{\Gamma(\chitu)> w}- e^w\pk{\Gamma(\eta^h )> w}\Biggr]\, dw }\\
&& \quad \leq   \limit{u}\sup_{\tau_u \in K_u}\int_{-M}^M
 e^w (1-   e^{-\frac{w^2}{2 \gtu^2}})\pk{\Gamma(\eta^h )> w} dw \\
 && \quad \quad  +\limit{u}\sup_{\tau_u \in K_u}\ABs{\int_{-M}^M
\Biggl[ e^{w-\frac{w^2}{2 \gtu^2}}\left(\pk{\Gamma(\chitu)> w}- \pk{\Gamma(\eta^h )> w}\right)\Biggr]\, dw }\\
&&\quad \leq   e^M \limit{u}\int_{-M}^M
\sup_{\tau_u \in K_u \COM{w\in [-M,M]}}\ABs{\pk{\Gamma(\chitu)> w}- \pk{\Gamma(\eta^h )> w} } dw =0.
\EQNY
}
Using \eqref{eqex} for $\delta \in (0,1/c)$, $|w|>M$ with $M$ sufficiently large and all $u$ large we have
$$ \sup_{\gE{\tau_u}\in K_u,t\in E} (1+h_{u,\tau_u}(t))\EE{\chitu({t})} \le \delta \abs{w}.$$
Moreover, in view of (\ref{eqco}) and (\ref{A31}) in {\bf C3} we have that  for $u$ sufficiently large
\BQNY
Var\Bigl( (1+h_{u,\tau_u}(t))\chitu({t})-(1+h_{u,\tau_u}(s))\chitu({s})\Bigr)
&\leq & \gtu^2\EE{\Bigl( Z_{u,\tau_u}({t})-Z_{u,\tau_u}({s})\Bigr)^2} \\
&  \leq&   \mathbb{Q}\sum_{i=1}^d|s_i-t_i|^a.
\EQNY Consequently, by  Piterbarg inequality
(see e.g., Theorem 8.1 in \cite{Pit96} \COM{and its extensions in Lemma 5.1 of \cite{KEP2016}, Theorem 3 of \cite{Pit16}})  we obtain for some $\ve \in (0,1)$, $\delta \in (0,1/c)$  with  $c$  given in {\bf F1}, and all $u$ large
\BQNY
\lefteqn{\int_{\abs{w}> M } e^{w-\frac{w^2}{2 \gtu^2}} \pk{\Gamma(\chitu)> w}\, dw}\\
&\le &
\int_{\abs{w}> M } e^{w}\pk{c\sup_{t\in E}  (1+h_{u,\tau_u}(t))(\chitu(t)- \EE{\chitu(t)})> w-c\sup_{t\in E, \tau_u \in K_u} (1+h_{u,\tau_u}(t))\EE{\chitu(t)} }\, dw\\
&\le & e^{-M}+ \int_{M}^\IF  e^{w}\Psi( (1- \ve)(1/c-\delta)w) \, dw \\
&=:&A(M) \to   0, \quad M\to  \IF.
\EQNY
\rdd{Moreover, by Borell-TIS inequality (see e.g., \cite{AdlerTaylor})
\BQNY  \int_{\abs{w}> M} e^{w}\pk{\Gamma( \eta^h)> w} dw
&\leq& \int_{\abs{w}> M} e^{w}\pk{c\sup_{t\in E}\ehd{\eta^h(t)}> w} dw\\
&\leq& e^{-M} +\int_{M}^\IF e^{w}\pk{\sqrt{2}c\sup_{t\in E}\eta(t)> w-c\sup_{t\in E}\left(Var(\eta(t))+h(t)\right)} dw\\
&\leq& e^{-M} + \int_{M}^\IF e^{w-\frac{(w-a)^2}{2\sup_{t\in E}Var(\sqrt{2}c\eta(t))}}  dw\\
&=:& B(M) \rw 0, \quad M\rw\IF,
\EQNY
with \rrd{$a=\sqrt{2}c\mathbb{E}\left\{\sup_{t\in E}\ehd{\eta(t)}\right\}-c\sup_{t\in E}\left(Var(\eta(t))+h(t)\right)<\IF$.}}
Hence \eqref{esB} follows from
\BQNY
 \sup_{\tau_u \in K_u} \abs{\mathcal{I}_{u,\tau_u}- \Hh(E) }&\le &
 \sup_{\tau_u \in K_u} \ABs{\int_{-M}^M
\Biggl[ e^{w-\frac{w^2}{2\gtu^2}}\pk{\Gamma(\chitu)> w}- e^w\pk{\Gamma( \eta^h )> w}\Biggr]\, dw }\\
&&
+ A(M)+ B(M)\\
&\to & A(M)+ B(M), \quad u\to \IF, \\
&\to &0, \quad M\to \IF,
\EQNY
establishing the proof.
\QED
\proofprop{ths} It follows from Remark \ref{remark} ii) that it suffices to prove (\ref{A1}), (\ref{A31}) and (\ref{eq3}). Without loss of generality, in the following derivation we assume that $c_i>0, 1\leq i\leq d$. By (\ref{A2}), we have
$$\theta_{u,\tau_u}(s,t)=\sum_{i=1}^d \frac{c_i\sigma_i^2(\qq_i(u)|s_i-t_i|)}{\sigma_i^2(\qq_i(u))}, \quad (s,t)\in E.$$
\ehd{By uniform convergence theorem (UCT) for regularly varying functions, see \cite{BI1989}}, (\ref{A1})  holds with $\eta$ defined in (\ref{A3}). Next we verify (\ref{A31}).  For  $0<\beta<\min(\min_{1\leq i\leq d}\alpha_{i,0}, \min_{d_2+1\leq i\leq d}\alpha_{i,\IF})$ we have
\BQNY
\sum_{i=1}^d \frac{c_i\sigma_i^2(\qq_i(u)|s_i-t_i|)}{\sigma_i^2(\qq_i(u))}= \sum_{i=1}^d c_i \frac{f_i(\qq_i(u)|s_i-t_i|)}{f_i(\qq_i(u))}|s_i-t_i|^{\beta/2},
\EQNY
with $f_i(t)=\frac{\sigma_i^2(t)}{t^{\beta/2}}, t>0$. Note that $f_i$ is regularly varying at $0$ with index $\alpha_{i,0}-\beta/2>0$ for $1\leq i\leq d$ and for $d_2+1\leq i\leq d$, $f_i$ is regularly varying at $\IF$ with index $\alpha_{i,\IF}-\beta/2>0$. By UCT for any $M>0$ we have
$$\lim_{u\rw\IF}\ehd{\max}_{i=1,\dots, d_1}\sup_{0<|s_i-t_i|\leq M}\left|\frac{f_i(\qq_i(u)|s_i-t_i|)}{f_i(\qq_i(u))}-|s_i-t_i|^{\alpha_{i,0}-\beta/2}\right|=0.$$
Using the fact that $f_i$ is bounded on compact intervals for $d_2+1\leq i\leq d$, again by UCT, for any $M>0$
$$\lim_{u\rw\IF}\ehd{\max}_{i=d_2+1,\dots, d}\sup_{0<|s_i-t_i|\leq M}\left|\frac{f_i(\qq_i(u)|s_i-t_i|)}{f_i(\qq_i(u))}-|s_i-t_i|^{\alpha_{i,\IF}-\beta/2}\right|=0.$$
Moreover, since $f_i$ is regularly varying at $0$ with index $\alpha_{i,0}-\beta>0$ and  $\varphi_i\in (0,\IF), d_1+1\leq i\leq d_2$, then  for  any $M>0$ and $u$ large enough
$$\ehd{\max}_{d_1+1\leq i\leq d_2}\sup_{0<|s_i-t_i|\leq M}\frac{f_i(\qq_i(u)|s_i-t_i|)}{f_i(\qq_i(u))}<\IF.$$
Thus we conclude that for $u$ large enough
$$\sum_{i=1}^d \frac{c_i\sigma_i^2(\qq_i(u)|s_i-t_i|)}{\sigma_i^2(\qq_i(u))}\leq \mathbb{Q}\sum_{i=1}^d|s_i-t_i|^{\beta/2}, \quad s,t\in E,$$
which confirms (\ref{A31}). We are now left to prove (\ref{eq3}). In light of (\ref{A2}) and {UCT}, we have
\BQNY&&\lim_{u\rw\IF} \sup_{t\in E\setminus\{0\}, {\tau_u} \in K_u}\ABs{ g_{u,\tau_u}^2Var(Z_{u,\tau_u}({t})-Z_{u,\tau_u}({0}))- 2 Var(\eta({t}))}\\
&& \quad \leq \lim_{u\rw\IF} \sup_{t\in E\setminus\{0\}, {\tau_u} \in K_u}\ABs{ \frac{g_{u,\tau_u}^2Var(Z_{u,\tau_u}({t})-Z_{u,\tau_u}({0}))}{2\theta_{u,\tau_u}(0,t)}-1}\ABs{2\theta_{u,\tau_u}(0,t)}\\
&& \ \ +\lim_{u\rw\IF} \sup_{t\in E, {\tau_u} \in K_u}\ABs{ 2\theta_{u,\tau_u}(0,t)- 2 Var(\eta({t}))}=0,
\EQNY
which implies that (\ref{eq3}) holds. This completes the proof.\QED
\prooftheo{cor} We check that {\bf C0-C3} hold. Clearly, \kd{{\bf C0} is satisfied  by the assumptions}. We observe that
$$\xitu(t)= \frac{\overline{\xi}_{u,\tau_u}(t)}{1+ h_{u,\tau_u}(t)}, \quad t\in E, \tau_u \in K_u,$$
with $$\overline{\xi}_{u,\tau_u}(t)=\frac{\xitu(t)}{\sigma_{u,\tau_u}(t)}, \ \ h_{u,\tau_u}(t)=\frac{1-\sigma_{u,\tau_u}(t)}{\sigma_{u,\tau_u}(t)},$$ which together with {\bf D1} immediately implies that {\bf C1} is valid. Let next for $u>0$
$$\theta_{u, \tau_u}(s, t)=\frac{g_{u,\tau_u}^2}{2}Var(\overline{\xi}_{u,\tau_u}(t)-\overline{\xi}_{u,\tau_u}(s)).$$
Direct calculations yield
\BQNY
\theta_{u, \tau_u}(s, t)=I_{1,u,\tau_u}(s,t)+I_{2,u,\tau_u}(s,t)+I_{3,u,\tau_u}(s,t),  \ \ s,t\in E,
\EQNY
where
$$I_{1,u,\tau_u}(s,t)=\frac{g_{u,\tau_u}^2}{2}\frac{Var(\xi_{u,\tau_u}(t)-\xi_{u,\tau_u}(s))}{\sigma_{u,\tau_u}^2(t)}, \ \ I_{2,u,\tau_u}(s,t)=\frac{g_{u,\tau_u}^2}{2}\frac{(\sigma_{u,\tau_u}(t)-\sigma_{u,\tau_u}(s))^2}
{\sigma_{u,\tau_u}^2(t)},$$
$$I_{3,u,\tau_u}(s,t)=g_{u,\tau_u}^2\frac{\sigma_{u,\tau_u}(t)-\sigma_{u,\tau_u}(s)}
{\sigma_{u,\tau_u}^2(t)\sigma_{u,\tau_u}(s)}\E{(\xi_{u,\tau_u}(s)-\xi_{u,\tau_u}(t))\xi_{u,\tau_u}(s)}.$$
It follows from {\bf D1} that
$$\lim_{u\rw\IF}\sup_{s,t\in E, \tau_u\in K_u}I_{2,u,\tau_u}(s,t)\leq \lim_{u\rw\IF}\sup_{s,t\in E, \tau_u\in K_u}g_{u,\tau_u}^2\frac{(\sigma_{u,\tau_u}(t)-1)^2+(1-\sigma_{u,\tau_u}(s))^2}
{\sigma_{u,\tau_u}^2(t)}\gE{=} 0.$$
Further, by {\bf D1,D2}
$$\lim_{u\rw\IF}\sup_{\tau_u\in K_u}\left|I_{1,u,\tau_u}(s,t)-Var(\eta(t)-\eta(s))\right| \gE{=} 0, \ \ s,t\in E
$$
and
$$\lim_{u\rw\IF}\sup_{\tau_u\in K_u}\left|I_{3,u,\tau_u}(s,t)\right|\leq \lim_{u\rw\IF}\sup_{\tau_u\in K_u}g_{u,\tau_u}^2\frac{|\sigma_{u,\tau_u}(t)-\sigma_{u,\tau_u}(s)|}
{\sigma_{u,\tau_u}^2(t)}\sqrt{Var(\xi_{u,\tau_u}(s)-\xi_{u,\tau_u}(t))} \gE{=0}, \ \ s,t\in E.$$
Thus we confirm that {\bf C2} holds.
\kd{Moreover,} by  {\bf D3} and the fact \KD{that}
$$(\sigma_{u,\tau_u}(t)-\sigma_{u,\tau_u}(s))^2\leq Var(\xi_{u,\tau_u}(t)-\xi_{u,\tau_u}(s))$$
\KD{we obtain}
$$\lim_{u\rw\IF}\sup_{ \tau_u\in K_u}\sup_{s\neq t, s,t\in E}\frac{\theta_{u, \tau_u}(s, t)}{||t-s||^\nu}\leq \mathbb{Q}\lim_{u\rw\IF}\sup_{ \tau_u\in K_u}\sup_{s\neq t, s,t\in E}\frac{g_{u,\tau_u}^2 Var(\xi_{u,\tau_u}(t)-\xi_{u,\tau_u}(s))}{||t-s||^\nu}<\IF.$$
Using again {\bf D1,D2} we obtain
$$\lim_{u\rw\IF}\sup_{t\in E,\tau_u\in K_u}\left|I_{1,u,\tau_u}(0,t)-Var(\eta(t))\right|=0,$$
$$  \lim_{u\rw\IF}\sup_{t\in E,\tau_u\in K_u}I_{2,u,\tau_u}(0,t)=0, \quad \lim_{u\rw\IF}\sup_{t\in E,\tau_u\in K_u}|I_{3,u,\tau_u}(0,t)|=0,$$ which imply
$$\lim_{u\rw\IF}\sup_{t\in E,\tau_u\in K_u}\left|\theta_{u,\tau_u}(0,t)-Var(\eta(t))\right|=0.$$
Hence {\bf C3} is satisfied with (\ref{eq3}) instead of (\ref{eq2}).
\kd{In view of  Remark \ref{remark} the proof is completed.} \QED

\prooftheo{pro} Recall that $F(A,B)= \inf_{s\in A, t\in B}\norm{s-t}$ with  $A, B$ \eM{two non-empty subsets of} $\mathbb{R}^d$ and $\norm{\cdot}$ the Euclidean norm. Clearly, for any $u$ positive
\BQNY
\pk{ \sup_{t\in \lambda_1+\mathcal{E}_1} X_u(t)> m_{\lambda_1}(u), \sup_{t\in\lambda_2 +\mathcal{E}_2} X_u(t)> m_{\lambda_2}(u)}
\leq \pk{ \sup_{s\in \lambda_1+\mathcal{E}_1, t\in \lambda_2 +\mathcal{E}_2}(X_u(s)+X_u(t))> 2m_{\lambda_1, \lambda_2}(u)},
\EQNY
where $m_{\lambda_1, \lambda_2}(u)=\min (m_{\lambda_1}(u), m_{\lambda_2}(u))$.
By (\ref{fu2}) and (\ref{delta}), we have that for $u$ sufficiently large and  $F(\lambda_1+\mathcal{E}_1, \lambda_2+\mathcal{E}_2)>S_1,$  with $S_1$ large enough,
\BQNY
2\delta\leq Var( X_u(s)+X_u(t))=4-2(1-r_u(s,t))\leq 4-\frac{2\mathcal{C}_1F^\beta(\lambda_1+\mathcal{E}_1, \lambda_2+\mathcal{E}_2)}{m^2(u)}.
\EQNY
Moreover, by (\ref{COR}) and the above inequality,
\BQNY\label{e1}
1-Cor(X_u(s)+X_u(t), X_u(s')+X_u(t'))&\leq& \frac{Var(X_u(s)+X_u(t)-X_u(s')-X_u(t'))}{2\sqrt{Var(X_u(s)+X_u(t))}\sqrt{Var(X_u(s')+X_u(t'))}}\\
&\leq& \delta^{-1}(1-r_u(s,s')+1-r_u(t,t'))\\
&\leq& \mathcal{C}_2\frac{\delta^{-1}d^{\alpha/2}}{m^2(u)}\sum_{i=1}^d(|s_i-s'_i|^{\alpha}+|t_i-t'_i|^{\alpha})
\EQNY
holds for $s,t,s',t'\in [0,1]^d$.
Let $X_u^*(s,t), s,t\in \mathbb{R}^d, u>0$ be  a family of centered Gaussian random fields with unit variance and correlation satisfying
$$r_u(s,t)=e^{-\frac{2\delta^{-1}d^{\alpha/2}\mathcal{C}_2}{m^2(u)}\sum_{i=1}^d(|s_i|^{\alpha}+|t_i|^{\alpha})}, \quad s,t\in \mathbb{R}^d
$$
and let further
$$m_{u,\lambda_1,\lambda_2,\mathcal{E}_1, \mathcal{E}_2}:=\frac{2m_{\lambda_1, \lambda_2}(u)}{\sqrt{4-\frac{2\mathcal{C}_1F^\beta(\lambda_1+\mathcal{E}_1, \lambda_2+\mathcal{E}_2)}{m^2(u)}}}, \quad I_{i_1,\dots, i_d}=\prod_{j=1}^d[i_j, i_j+1].$$
For all $u$ large we have
\BQN
\lefteqn{ \pk{ \sup_{s\in \lambda_1+\mathcal{E}_1, t\in \lambda_2 +\mathcal{E}_2}( X_u(s)+X_u(t))> 2m_{\lambda_1, \lambda_2}(u)}}\nonumber\\
& \leq &
\pk{ \sup_{s\in \lambda_1+\mathcal{E}_1, t\in \lambda_2 +\mathcal{E}_2}  \overline{X_u(s)+X_u(t)}> m_{u,\lambda_1,\lambda_2,\mathcal{E}_1, \mathcal{E}_2} }\nonumber\\
&\leq& \pk{ \sup_{s\in \lambda_1+[0,S_2]^d, t\in \lambda_2 +[0,S_2]^d}  \overline{X_u(s)+X_u(t)}> m_{u,\lambda_1,\lambda_2,\mathcal{E}_1, \mathcal{E}_2} }\nonumber\\
&\leq& \sum_{i_1, i_2,\dots, i_d, i_1', i_2',\dots, i_d'=0}^{[S_2]}\pk{ \sup_{s\in \lambda_1+I_{i_1,\dots, i_d}, t\in \lambda_2 +I_{i_1',\dots, i_d'}}  \overline{X_u(s)+X_u(t)}> m_{u,\lambda_1,\lambda_2,\mathcal{E}_1, \mathcal{E}_2} }\nonumber\\
&\leq& \sum_{i_1, i_2,\dots, i_d, i_1', i_2',\dots, i_d'=0}^{[S_2]}\pk{ \sup_{s\in \lambda_1+I_{i_1,\dots, i_d}, t\in \lambda_2 +I_{i_1',\dots, i_d'}}  X_u^*(s,t)> m_{u,\lambda_1,\lambda_2,\mathcal{E}_1, \mathcal{E}_2} }\nonumber\\
&=& (S_2+1)^{2d} \pk{ \sup_{s,t\in [0,1]^d}  X_u^*(s,t)> m_{u,\lambda_1,\lambda_2,\mathcal{E}_1, \mathcal{E}_2} },\label{new0}
\EQN
where we used Slepian inequality (see, e.g., \cite{AdlerTaylor,AZI}) to derive (\ref{new0}).
Hence in order to complete the proof, we need to apply Proposition \ref{ths} to the family of Gaussian random fields $\{X_u^*(s,t), (s,t)\in [0,1]^{2d}\}$. Let
$$  K_u=\{(\lambda_1, \lambda_2), \lambda_i+\mathcal{E}_i\subset E_u, i=1,2\}. \quad $$
Note that
$$
\lim_{u\rw\IF}\sup_{(\lambda_1, \lambda_2)\in K_u}\sup_{(s,t)\neq (s',t'), (s,t), (s',t')\in [0,1]^{2d}}\left| \frac{(m_{u,\lambda_1,\lambda_2,\mathcal{E}_1, \mathcal{E}_2} ) ^2 Var\left(X_u^*(s,t)-X_u^*(s',t')\right)}
{2\sum_{i=1}^d 2\delta^{-1}d^{\alpha/2}\mathcal{C}_2(\sum_{i=1}^d|s_i-s_i'|^{\alpha}+\sum_{i=1}^d|t_i-t_i'|^{\alpha})}-1\right| = 0.
$$
Since \eM{conditions {\bf C0-C1} are clearly satisfied,} then  Proposition \ref{ths} implies
\BQNY
\lim_{u\rw\IF}\sup_{(\lambda_1,\lambda_2)\in K_u}\left| \frac 1 { \Psi\left(m_{u,\lambda_1,\lambda_2,\mathcal{E}_1, \mathcal{E}_2}\right)}\pk{ \sup_{s,t\in [0,1]^{2d}} X_u^*(s,t)> m_{u,\lambda_1,\lambda_2,\mathcal{E}_1, \mathcal{E}_2}}-\mathcal{H}_\eta([0,1]^{2d})\right|=0,
\EQNY
where
$$\eta(s,t)=\sum_{i=1}^{d} \sqrt{2\delta^{-1}d^{\alpha/2}\mathcal{C}_2}B_{\alpha}^{(i)}(s_i)+\sum_{i=d+1}^{2d}\sqrt{2\delta^{-1}d^{\alpha/2}\mathcal{C}_2}B_{\alpha}^{(i)}(t_{i-d}), $$
with $B_{\alpha}^{(i)}, 1\leq i\leq 2d$ independent fBm's with index $\alpha$.
Thus we  establish the claim for $F(\lambda_1+\mathcal{E}_1, \lambda_2+\mathcal{E}_2)>S_1$. For $F(\lambda_1+\mathcal{E}_1, \lambda_2+\mathcal{E}_2)\leq S_1$, we have
\BQNY
\pk{ \sup_{s\in \lambda_1+\mathcal{E}_1} X_u(s)> m_{\lambda_1}(u), \sup_{t\in\lambda_2 +\mathcal{E}_2} X_u(t)> m_{\lambda_2}(u)}\leq
\pk{ \sup_{t\in \lambda_1+[-S_1,S_2+S_1]^d} X_u(t)> m_{\lambda_1,\lambda_2}(u)}.
\EQNY
By (\ref{COR}) and Slepian inequality
\BQNY
&&\pk{ \sup_{s\in \lambda_1+[-S_1,S_2+S_1]^d} X_u(s)> m_{\lambda_1,\lambda_2}(u)}\\
&& \ \ \leq (S_2+2S_1+1)^d \pk{ \sup_{s\in [0,1]^d} X_u^*(\delta^{1/\alpha}s, 0,\dots, 0)> m_{\lambda_1,\lambda_2}(u)}\\
&& \ \ \sim (S_2+2S_1+1)^d \mathcal{H}_{\lambda}([0,1]^d)\Psi(m_{\lambda_1,\lambda_2}(u)), \ \ u\rw\IF,
\EQNY
with $\lambda(s)=\sqrt{\delta}\eta(s, 0 \ldot  0)$.
This completes the proof.
\QED
\proofkorr{CC0} Let
$\beta=\frac{1}{2}\min_{i=1,\dots, 2d}\min(\alpha_{i,0},
\alpha_{i,\IF}, 2)$ and  $f_i(t)=\kappa_i(t)/t^\beta$.
 Clearly, $f_i$'s \ehd{are}  regularly varying at $0$ with index $\alpha_{i,0}-\beta>0$
 and regularly varying at $\IF$ with {index}  $\alpha_{i,\IF}-\beta>0$.  \ehd{With this notation we have}
\bqn{\frac{\kappa_i(\qq_i(u)|s_i-t_i|)}{\kappa_i(\qq_i(u))}=\frac{f_i(\qq_i(u)|s_i-t_i|)}{f_i(\qq_i(u))}|s_i-t_i|^\beta, \quad s_i\neq t_i, i=1,\dots, 2d.}
 Next we focus on $\frac{f_i(\qq_i(u)|s_i-t_i|)}{f_i(\qq_i(u))}$. We consider the upper bound and lower bound respectively.\\
 {\it \underline{Lower bound}}. For $\varphi_i=0$ we define  $g_i(t)=1/ f_i(1/t)$. Then $g_i$ is both regularly varying at $0$ with index $\alpha_{i,\IF}-\beta>0$ and \ehd{regularly varying} at $\IF$ with index $\alpha_{i,0}-\beta>0$. By the assumption on $\kappa_i$'s, further  $g_i$ is bounded over any compact interval and by UCT
 $$\lim_{u\rw\IF}\sup_{|s_i-t_i|\geq 1}\left|\frac{g_i(\frac{1}{\qq_i(u)|s_i-t_i|})}{g_i(\frac{1}{\qq_i(u)})}-\left(\frac{1}{|s_i-t_i|}\right)^{\alpha_{i,0}-\beta}\right|=0$$
 implying that for $u$ large enough
 $$\frac{g_i(\frac{1}{\qq_i(u)|s_i-t_i|})}{g_i(\frac{1}{\qq_i(u)})}\leq 2, \quad \frac{1}{|s_i-t_i|}\leq 1. $$
 Consequently, \kd{for $u$ sufficiently large}
$$
	 \frac{f_i(\qq_i(u)|s_i-t_i|)}{f_i(\qq_i(u))}=\frac{g_i(\frac{1}{\qq_i(u)})}{g_i(\frac{1}{\qq_i(u)|s_i-t_i|})}\geq \frac{1}{2}, \quad |s_i-t_i|\geq 1.
	 $$
Next, if  $\varphi_i\in (0,\IF)$, then by the fact that $\lim_{t\rw\IF}f_i(t)=\IF$, there exists $S_1>0$ and $M_i'$ such that for $u$ sufficiently large
$$\frac{f_i(\qq_i(u)|s_i-t_i|)}{f_i(\qq_i(u))}>M_i', \quad |s_i-t_i|>S_1.$$
For $\varphi=\IF$, Potter's theorem (see e.g., \cite{BI1989}[Theorem 1.5.6]) implies that for any $0<\epsilon<\alpha_{i,\IF}-\beta$ there exists $M_i''>0$ and $S_1'>1$ such that for $u$ sufficiently large
$$\frac{f_i(\qq_i(u)|s_i-t_i|)}{f_i(\qq_i(u))}>M_i''|s_i-t_i|^{\alpha_{i,\IF}-\beta-\epsilon}\geq M_1'', \quad |s_i-t_i|>S_1'.$$
Consequently, there exists $S>1$ and $M>0$ such that  for $u$ sufficiently large
$$\frac{\kappa_i(\qq_i(u)|s_i-t_i|)}{\kappa_i(\qq_i(u))}\geq M|s_i-t_i|^\beta, \quad  |s_i-t_i|>S, i=1,\dots, d.$$
Further, for $u$ large enough
\BQN\label{low11}
g_u(s,t)\geq d^{-\frac{\beta}{2}}M\norm{s-t}^\beta, \quad  \norm{s-t}>\sqrt{d}S.
\EQN
{\it \underline{Upper bound}}. If $\varphi_i \in \{ 0, \IF\}$, then using again UCT we have that
$$\ehd{\sup_{|s_i-t_i|\leq 1}}\frac{f_i(\qq_i(u)|s_i-t_i|)}{f_i(\qq_i(u))}\leq \ehd{C}$$
is valid for all $u$ large enough  and some constant $C$. \ehd{Further, since $f_i$ \kd{is}  locally bounded}, then the above holds
also if $\varphi_i\in (0,\IF)$. This implies that for some $M'>0$
\BQN
\widetilde{g}_u(s,t)\leq M'\sum_{i=1}^d|s_i-t_i|^\beta\leq d M'\norm{s-t}^\beta, \quad s-t\in [-1,1]^d,
\EQN
which combined with (\ref{low11}) and Theorem \ref{pro} establishes the claim. \QED
\proofkorr{CC1} The claim follows straightforwardly using the arguments of  Corollary \ref{CC0} for the case $\varphi_i=0$. \QED

\prooftheo{thm} Without loss of generality, we assume that $a_i=-\IF$, $b_i=\IF$ for $d_1+1\leq i\leq d_2$.
Set in the following
$$I_k=\prod_{i=1}^{d_1}[k_iS, (k_i+1)S], \quad k=(k_1,\dots, k_{d_1}),$$
$$J_l=\prod_{i={d_1+1}}^{d_2}[l_iS, (l_i+1)S]\times \prod_{i={d_2+1}}^{d}[l_iT, (l_i+1)T], \quad l=(l_{d_1+1}, \dots, l_{d}),$$ $$J^*=\prod_{i={d_1+1}}^{d_2}[-S,S]\times \prod_{i={d_2+1}}^{d}[-T,T], \widetilde{J}=\prod_{i={d_1+1}}^{d_2}[-S,S]\times\{0\}, \quad 0\in \mathbb{R}^{d-d_2}.$$
Further, define
$$ I_k^*=I_k\times J^*\times \HEH{E} , \quad
\widetilde{I}_k=I_k\times \widetilde{J}\times \HEH{E}, \quad {I_{k,l}=I_k\times J_l\times E},
$$
$$K_u^{\pm}=\Biggl\{k, \frac{a_i(u)}{S}\mp 1\leq k_i\leq \frac{b_i(u)}{S}\pm 1, 1\leq i\leq d_1 \Biggr\},
$$
$$L_u=\Biggl\{l, \frac{a_i(u)}{S}- 1\leq l_i\leq \frac{b_i(u)}{S}+ 1,  d_1+1\leq i\leq d_2,
 \quad \frac{a_i(u)}{T}- 1\leq l_i\leq \frac{b_i(u)}{T}+ 1,  d_2+1\leq i\ehd{\leq d}, {J_l\nsubseteqq J^*} \Biggr\}.$$
For some $\ve\in (-1,1)$ and $u>0$ set
$$\Theta_{\epsilon}(u):=\prod_{i=1}^{d_1}\int_{y_{i,1}}^{y_{i,2}}e^{-(1-\epsilon)|s|^{\beta_i}}ds \prod_{i=1}^{d_1}\left(\frac{g_i(u)}{m^2(u)}\right)^{1/\beta_i}\Psi(m(u)).$$
Observe that
$$X_u(s,t)=\frac{\sigma_u(s,t)\overline{X}_u(s,t)}{\sigma_u(0,0)}, \quad \frac{\sigma_u(0,0)}{\sigma_u(s,t)}=\frac{\sigma_u(0,0)}{\sigma_u(s,0)}  \frac{\sigma_u(s,0)}{\sigma_u(s,t)}.$$
Using (\ref{var2}) and (\ref{new1}), there exists $e_{u,1}(s)$ and $e_{u,2}(s,t)$ such that as $u\rw\IF$
$$\sup_{s\in \prod_{i=1}^d[a_i(u), b_i(u)]}|e_{u,1}(s)|=o(1), \quad \sup_{(s,t)\in E_u}|e_{u,2}(s,t)|=o(1),$$
\rdd{and
\BQNY
\frac{\sigma_u(0,0)}{\sigma_u(s,0)}&=&1+(1+e_{u,1}(s))\sum_{i=1}^d\frac{|s_i|^{\beta_i}}{g_i(u)}, \quad  s\in \prod_{i=1}^d[a_i(u), b_i(u)], \\ \frac{\sigma_u(s,0)}{\sigma_u(s,t)}&=&1+(1+e_{u,2}(s,t))
\sum_{i=d+1}^{d+n}\frac{|t_i|^{\beta_i}}{g_i(u)}, \quad (s,t)\in E_u.
\EQNY
Note that by {\bf F2} for $\Gamma^*$
$$\Gamma_{E_u}(X_u(s,t))=\sup_{s\in \prod_{i=1}^d[a_i(u), b_i(u)]}\Gamma^*(X_u(s,t))=\sup_{s\in \prod_{i=1}^d[a_i(u), b_i(u)]}\sigma_u(s,0)\Gamma^*\left(\overline{X}_u(s,t)\frac{\sigma_u(s,t)}{\sigma_u(s,0)}\right).$$}
Thus, by  {\bf F2} for $\Gamma^*$, and the property of $\sup$ functional we have that  for  $0<\epsilon<1/2$ and $u$ sufficiently large
\BQN\label{dec}
\pk{\Gamma_{E_u}(X_u^{+\epsilon})>m(u)}\leq \pk{\Gamma_{E_u}(X_u)>m(u)}\leq \pk{\Gamma_{E_u}(X_u^{-\epsilon,y})>m(u)},
\EQN
where for $(s,t)\in E_u$
\BQNY
X_u^{-\epsilon,y}(s,t)&=&\frac{\overline{X}_u(s,t)}{(1+\sum_{i=1}^{d_1}(1-\epsilon)\frac{|s_i|^{\beta_i}}
{g_i(u)})(1+\sum_{i=d_1+1}^{d_2}(1-\epsilon)\frac{|s_i|^{\beta_i}}{g_i(u)}+\sum_{i=d_2+1}^{d}y\frac{|s_i|^{\beta_i}}
{m^2(u)})}\\
&& \ \ \times \frac{1}{(1+(1+e_{u,2}(s,t))\sum_{i=d+1}^{d+n}\frac{|t_i|^{\beta_i}}{g_i(u)})},
\EQNY
and
$$X_{u}^{+\epsilon}(s,t)=\frac{\overline{X}_u(s,t)}{(1+\sum_{i=1}^{d_1}(1+\epsilon)\frac{|s_i|^{\beta_i}}
{g_i(u)})(1+\sum_{i=d_1+1}^{d}(1+\epsilon)\frac{|s_i|^{\beta_i}}{g_i(u)})(1+(1+e_{u,2}(s,t))\sum_{i=d+1}^{d+n}\frac{|t_i|^{\beta_i}}{g_i(u)})}.$$
{\it \underline{Upper bound}}.
By the  property of $\sup$ functional, we have that
\BQN\label{add}
\pk{\Gamma_{E_u}(X_u^{-\epsilon,y})>m(u)}&\leq& \sum_{k\in K_u^+}
\pk{\Gamma_{I_k^*}(X_u^{-\epsilon,y})>m(u)}
+\sum_{(k,l)\in K_u^+\times L_u }\pk{\Gamma_{{I_{k,l}}}(X_u^{-\epsilon,y})>m(u)}\nonumber\\
&\leq & \sum_{k\in K_u^+}\pk{\Gamma_{I_0^*}(\xi_{u,k})>m_{u,k}}+\sum_{(k,l)\in K_u^
+\times L_u }\pk{\Gamma_{{I_{0,0}}}(\xi_{u,k,l})>m_{u,k,l}},
\EQN
where
\BQNY\xi_{u,k}(s,t)&=&\frac{\overline{X}_u(s+kS,t)}{(1+\sum_{i=d_1+1}^{d_2}(1-\epsilon){\frac{|s_i|^{\beta_i}}{g_i(u)}}+\sum_{i=d_2+1}^{d}y
\frac{|s_i|^{\beta_i}}
{m^2(u)})(1+(1+e_{u,2}(s,t))\sum_{i=d+1}^{d+n}\frac{|t_i|^{\beta_i}}{g_i(u)})}, \quad (s,t)\in I_0^*,\\
\xi_{u,k,l}(s,t)&=&\frac{\overline{X}_u(s+(k,l)(S,T),t)}{1+(1+e_{u,2}(s,t))\sum_{i=d+1}^{d+n}\frac{|t_i|^{\beta_i}}{g_i(u)}},  \quad (s,t)\in {I_{0,0}},\\
m_{u,k}&=&m(u)\left(1+\sum_{i=1}^{d_1}(1-\epsilon)\frac{|k_i^*S|^{\beta_i}}
{g_i(u)}\right),\\
 m_{u,k,l}&=&m(u)\left(1+\sum_{i=1}^{d_1}(1-\epsilon)\frac{|k_i^*S|^{\beta_i}}
{g_i(u)}+\sum_{i=d_1+1}^{d_2}(1-2\epsilon)\frac{|l_i^*S|^{\beta_i}}{g_i(u)}+\sum_{i=d_2+1}^{d}y/2\frac{|l_i^*S|^{\beta_i}}
{m^2(u)}\right),
\EQNY
with $kS=(k_1S, \dots, k_{d_1}S, 0, \dots ,0)\in \mathbb{R}^d$ and
$$(k,l)(S,T)=(k_1S, \dots, k_{d_1}S, l_{d_1+1}S, \dots ,l_{d_2}S, l_{d_2+1}T, l_{d}T)\in \mathbb{R}^d,$$
$$ k_i^*=\min(|k_i|, |k_i+1|), \quad  1\leq i\leq d_1,  l_i^*=\min(|l_i|, |l_i+1|), \quad d_1+1\leq i\leq d.$$
In order to apply Proposition \ref{ths}, by (\ref{cor2}), set
$$\theta_{u,k}(s,t,s',t')=\sum_{i=1}^d{\frac{c_i\sigma_i^2(\qq_i(u)|s_i-s_i'|)}{\sigma^2_i(\qq_i(u))}}+\sum_{i=d+1}^{d+n}
\frac{c_i\sigma_i^2(\qq_i(u)|t_i-t_i'|)}{\sigma_i^2(\qq_i(u))}, \quad  (s,t),  (s',t')\in I_0^*,$$
$$ h_{u,k}(s,t)=\left(\sum_{i=d_1+1}^{d_2}(1-\epsilon)\frac{|s_i|^{\beta_i}}{g_i(u)}+\sum_{i=d_2+1}^{d}y\frac{|s_i|^{\beta_i}}
{m^2(u)}+\sum_{i=d+1}^{d+n}\frac{|t_i|^{\beta_i}}{g_i(u)}\right)(1+o(1)),  \quad (s,t)\in I_0^*,$$
$$g_{u,k}=m_{u,k}, \quad K_u=K_u^+, \quad E=I_0^*.$$
First we note that condition {\bf C0} holds straightforwardly. One can easily check that {\bf C1} holds with
\BQN\label{h}
h_\epsilon(s,t)=\sum_{i=d_1+1}^{d_2}(1-\epsilon)\gamma_i|s_i|^{\beta_i}+\sum_{i=d_2+1}^{d}y|s_i|^{\beta_i}
+\sum_{i=d+1}^{d+n}\gamma_i|t_i|^{\beta_i}, \quad (s,t)\in I_0^*.
\EQN
 Thus \rdd{in view of {\bf A1-A2} } and  by Proposition \ref{ths}, we have
\BQN\label{uni}
\lim_{u\rw\IF}\sup_{k\in K_u^+}\left|\frac{\pk{\Gamma_{I_0^*}(\xi_{u,{k}})>m_{u,k}}}{\Psi(m_{u,k})}-\mathcal{H}_{V_\varphi,  h_\epsilon}^{\Gamma}(I_0^*)\right|=0,
\EQN
with $h_\epsilon$ defined in (\ref{h}) and {$V_\varphi(s,t)=\sum_{i=1}^dV_{\varphi_i}(s_i)+\sum_{i=1}^n V_{\varphi_{d+i}}(t_i)$} with $V_{\varphi_i}$ defined in (\ref{vf}).
Similarly, we have
\BQN
\lim_{u\rw\IF}\sup_{(k,l)\in K_u^+\times L_u}\left|\frac{\pk{\Gamma_{{I_{0,0}}}(\xi_{u,k,l})>m_{u,k,l}}}{\Psi(m_{u,k,l})}-\mathcal{H}_{V_\varphi, \widetilde{h}}^{\Gamma}({I_{0,0}})\right|=0,
\EQN
with $\widetilde{h}(s,t)=\sum_{i=1}^{n}\gamma_{i+d}|t_i|^{\beta_{i+d}}$. Further, as $u\rw\IF$
\BQN\label{upper1}
\sum_{k\in K_u^+}\pk{\Gamma_{I_0^*}(\xi_{u,k})>m_k(u)}&\sim& \mathcal{H}_{V_\varphi,  h_\epsilon}^{\Gamma}(I_0^*)\sum_{k\in K_u^+}\Psi(m_{u,k})\nonumber\\
&\sim& \mathcal{H}_{V_\varphi,  h_\epsilon}^{\Gamma}(I_0^*)\Psi(m(u))\sum_{k\in K_u^+}e^{-\sum_{i=1}^{d_1}(1-\epsilon)m^2(u)\frac{|k_i^*S|^{\beta_i}}
{g_i(u)}}\nonumber\\
&\sim& S^{-d_1}\mathcal{H}_{V_\varphi,  h_\epsilon}^{\Gamma}(I_0^*) \Theta_{\epsilon}(u)
\EQN
and
\BQN\label{upper2}
&&\sum_{(k,l)\in K_u^+\times L_u }\pk{\Gamma_{I_{0,0}}(\xi_{u,k,l})>m_{u,k,l}}\nonumber\\
&& \ \ \sim \mathcal{H}_{V_\varphi, \widetilde{h}}^{\Gamma}(I_{0,0})\sum_{(k,l)\in K_u^+\times L_u }\Psi(m_{u,k,l})\nonumber\\
&& \ \ \leq\mathcal{H}_{V_\varphi, \widetilde{h}}^{\Gamma}(I_{0,0})\sum_{k\in K_u^+}\Psi(m_{u,k})\sum_{l\in L_u } e^{-m^2(u)(\sum_{i=d_1+1}^{d_2}(1-2\epsilon)\frac{|l_i^*S|^{\beta_i}}{g_i(u)}+\sum_{i=d_2+1}^{d}y/2\frac{|l_i^*T|^{\beta_i}}{m^2(u)})}(1+o(1))\nonumber\\
&& \ \ \leq\mathcal{H}_{V_\varphi, \widetilde{h}}^{\Gamma}(I_{0,0})\sum_{k\in K_u^+}\Psi(m_{u,k})\sum_{l\in L_u } e^{-\sum_{i=d_1+1}^{d_2}(1-2\epsilon)\gamma_i|l_i^*S|^{\beta_i}-\sum_{i=d_2+1}^{d}y/2|l_i^*T|^{\beta_i}}(1+o(1))\nonumber\\
&& \ \ \leq S^{-d_1}\mathcal{H}_{V_\varphi, \widetilde{h}}^{\Gamma}(I_{0,0}) \left(\sum_{i=d_1+1}^{d_2}e^{-\mathbb{Q}S^{\beta_i}}+\sum_{i=d_2+1}^{d}e^{-y\mathbb{Q}T^{\beta_i}}\right){\Theta_{\epsilon}(u)}(1+o(1)).
\EQN
{\it \underline{Lower bound}}.
By the property of $\sup$ functional and Bonferroni inequality, we obtain
\BQN\label{dec1}
\pk{\Gamma_{E_u}(X_u^{+\epsilon})>m(u)}&\geq& \sum_{k\in K_u^-}\pk{\Gamma_{\widetilde{I}_k}(X_u^{+\epsilon})>m(u)}\nonumber\\
&& \ \ -\sum_{k, q\in K_u^-, k\neq q}
\pk{\Gamma_{\widetilde{I}_k}(X_u^{+\epsilon})>m(u), \Gamma_{\widetilde{I}_q}(X_u^{+\epsilon})>m(u)}.
\EQN
Similarly as (\ref{upper1}), we have
\BQN\label{low1}
&&\sum_{k\in K_u^-}\pk{\Gamma_{\widetilde{I}_k}(X_u^{+\epsilon})>m(u)}\sim S^{-d_1}\mathcal{H}_{V_\varphi, h^*_\epsilon}^{\Gamma}(\widetilde{I}_0)\Theta_{-\epsilon}(u) ,
\EQN
with $h^*_\epsilon(s,t)=\sum_{i=d_1+1}^{d_2}(1+\epsilon)\gamma_i|s_i|^{\beta_i}+\sum_{i=1}^{n}\gamma_{i+d}|t_i|^{\beta_{i+d}}, (s,t)\in \widetilde{I}_0$.
Finally, we focus on the \ehd{double-sum} term. It follows from {\bf F1}, that
\BQNY
&&\sum_{k, q\in K_u^-, k\neq q}
\pk{\Gamma_{\widetilde{I}_k}(X_u^{+\epsilon})>m(u), \Gamma_{\widetilde{I}_k}(X_u^{+\epsilon})>m(u)}\\
&&\leq \sum_{k, q\in K_u^-, k\neq q}
\pk{\sup_{(s,t)\in\widetilde{I}_k} X_u^{+\epsilon}(s,t)>m(u), \sup_{(s,t)\in\widetilde{I}_q}X_u^{+\epsilon}(s,t)>m(u)}\\
&&\leq \sum_{k, q\in K_u^-, k\neq q}
\pk{\sup_{(s,t)\in\widetilde{I}_k} \overline{X}_u(s,t)>m_{u,k}, \sup_{(s,t)\in\widetilde{I}_q}\overline{X}_{u}(s,t)>m_{u,q}}.
\EQNY
Let for $u>0$
$$\mathcal{T}_1=\{(k,q), k, q\in K_u^-, k\neq q,  \widetilde{I}_k\kd{\cap} \widetilde{I}_q \neq \emptyset\},
\quad \mathcal{T}_2=\{(k,q), k, q\in K_u^-,  \widetilde{I}_k\kd{\cap} \widetilde{I}_q = \emptyset\}.$$
 Without loss of generality, we assume that $q_1=k_1+1, S>1$. Then $\widetilde{I}_k=\widetilde{I}_k'\bigcup \widetilde{I}_k''$ with $$\widetilde{I}_k'=[k_1S, (k_1+1)S-\sqrt{S}]\times \prod_{i=2}^{d_1}[k_iS,(k_i+1)S]\times\widetilde{J}\times E,$$
$$ \widetilde{I}_k''=[ (k_1+1)S-\sqrt{S}, (k_1+1)S]\times \prod_{i=2}^{d_1}[k_iS,(k_i+1)S]\times\widetilde{J}\times E.$$
Consequently,
\BQNY
&&\pk{\sup_{(s,t)\in\widetilde{I}_k} \overline{X}_u(s,t)>m_{u,k}, \sup_{(s,t)\in\widetilde{I}_q}\overline{X}_{u}(s,t)>m_{u,q}}\nonumber\\
&& \ \ \leq \pk{\sup_{(s,t)\in\widetilde{I}_k'} \overline{X}_u(s,t)>m_{u,k}, \sup_{(s,t)\in\widetilde{I}_q}\overline{X}_{u}(s,t)>m_{u,q}}
+\pk{\sup_{(s,t)\in\widetilde{I}_k''} \overline{X}_u(s,t)>m_{u,k}}.
\EQNY
Similarly as in (\ref{uni}), we have
\BQNY
\lim_{u\rw\IF}\sup_{k\in K_u^-}\left|\frac{\pk{\sup_{(s,t)\in\widetilde{I}_k''} \overline{X}_u(s,t)>m_{u,k}}}{\Psi(m_{u,k})}-\mathcal{H}_{V_\varphi, h^*_\epsilon}^{\sup}(\widehat{I}_0)\right|=0,
\EQNY
with $\widehat{I}_0=[0,\sqrt{S}]\times[0,S]^{d_1-1}\times\widetilde{J}\times E$.

 Let $\beta=\min(\min_{i=1}^{d+n}\alpha_{i,0}, \min_{i=1}^{d+n} \alpha_{i,\IF})$. By  (\ref{cor2}) and \rd{Corollary} \ref{CC0}, there exists $\mathcal{C}>0$ and $\mathcal{C}_1>0$ such that
 \BQNY
 &&\pk{\sup_{(s,t)\in\widetilde{I}_k'} \overline{X}_u(s,t)>m_{u,k}, \sup_{(s,t)\in\widetilde{I}_q}\overline{X}_{u}(s,t)>m_{u,q}}\\
 && \ \ \leq \mathcal{C}(S+|E|+1)^{2(d_2+n)}e^{-\mathcal{C}_1S^{\beta/2}}\Psi(m_{u,k,q}^*)
 \EQNY
 and for $(k,q)\in \mathcal{T}_2$
 \BQNY
 &&\pk{\sup_{(s,t)\in\widetilde{I}_k} \overline{X}_u(s,t)>m_{u,k}, \sup_{(s,t)\in\widetilde{I}_q}\overline{X}_{u}(s,t)>m_{u,q}}\\
 && \ \ \leq \mathcal{C}(S+|E|+1)^{2(d_2+n)}e^{-\mathcal{C}_1F^\beta(\widetilde{I}_k, \widetilde{I}_q)}\Psi(m_{u,k,q}^*),
 \EQNY
with $m_{u,k,q}^*=\min (m_{u,k}, m_{u,q})$. Since each $\widetilde{I}_k$ has at most $3^{d_1}$ neighbours, then for $S$ and $u$ sufficiently large
\BQN\label{low2}
&&\sum_{(k,q)\in \mathcal{T}_1}\pk{\sup_{(s,t)\in\widetilde{I}_k} \overline{X}_u(s,t)>m_{u,k}, \sup_{(s,t)\in\widetilde{I}_q}\overline{X}_{u}(s,t)>m_{u,q}}\nonumber\\
&& \ \ \leq 3^d\sum_{k\in K_u^-}\mathcal{H}_{V_\varphi, h_\epsilon^*}^{\sup}(\widehat{I}_0)\Psi(m_{u,k})+\sum_{(k,q)\in \mathcal{T}_1}\mathcal{C}(S+|E|+1)^{2(d_2+n)}e^{-\mathcal{C}_1
S^{\beta/2}}\Psi(m_{u,k,q}^*)\nonumber\\
&& \ \ \leq \mathbb{Q}\sum_{k\in K_u^-}\left(\mathcal{H}_{V_\varphi, h_\epsilon^*}^{\sup}(\widehat{I}_0)+e^{-\frac{\mathcal{C}_1
S^{\beta/2}}{2}}\right)\Psi(m_{u,k})\nonumber\\
&& \ \ \leq \mathbb{Q}S^{-d_1}\left(\mathcal{H}_{V_\varphi, h_\epsilon^*}^{\sup}(\widehat{I}_0)+e^{-\frac{\mathcal{C}_1
S^{\beta/2}}{2}}\right)\Theta_{\epsilon}(u).
\EQN
Moreover, for all $u$ large
\BQN\label{low3}
&&\sum_{(k,q)\in \mathcal{T}_2}\pk{\sup_{(s,t)\in\widetilde{I}_k} \overline{X}_u(s,t)>m_{u,k}, \sup_{(s,t)\in\widetilde{I}_q}\overline{X}_{u}(s,t)>m_{u,q}}\nonumber\\
&& \ \ \leq \sum_{(k,q)\in \mathcal{T}_2} \mathcal{C}(S+|E|+1)^{2(d_2+n)}e^{-\mathcal{C}_1 F^\beta(\widetilde{I}_k, \widetilde{I}_q)}\Psi(m_{u,k,q})\nonumber\\
&& \ \ \leq \sum_{k\in K_u^-} \Psi(m_{u,k}) \mathbb{Q}S^{\mathbb{Q}_1}\sum_{q\neq 0}e^{-\mathcal{C}_1(S^2\sum_{i=1}^{d_1}q_i^2)^{\beta/2}}\notag\\
&& \ \ \leq \mathbb{Q}S^{\mathbb{Q}_1}e^{-\mathbb{Q}_2S^{{\beta}}}\Theta_{\epsilon}(u).
\EQN
Inserting (\ref{add}--\ref{low3}) into (\ref{dec}) and dividing each term by $\Theta_{0}(u)$,
we have, with $\epsilon\rw 0$
\BQN\label{main}
&&S^{-d_1}\mathcal{H}_{V_\varphi, h_0^*}^{\Gamma}(\widetilde{I}_0)-\mathbb{Q} S^{-d_1}\left(\mathcal{H}_{V_\varphi, h_0^*}^{\sup}(\widehat{I}_0)+e^{-\frac{\mathcal{C}_1
S^{\beta/2}}{2}}\right)-\mathbb{Q}S^{\mathbb{Q}_1}e^{-\mathbb{Q}_2S^{\beta}}\nonumber\\
&& \ \ \leq \liminf_{u\rw\IF}\frac{\pk{\Gamma_{E_u}(X_u)>m(u)}}{\Theta_0(u)} \notag \\
&& \ \ \leq\lim_{T\rw 0}\lim_{y\rw\IF}\limsup_{u\rw\IF}\frac{\pk{\Gamma_{E_u}(X_u)>m(u)}}{\Theta_0(u)}\nonumber\\
&& \ \ \leq \lim_{T\rw 0}S^{-d_1}\mathcal{H}_{V_\varphi, h_0}^{\Gamma}(I_0^*)+\lim_{T\rw 0}\lim_{y\rw\IF}S^{-d_1}\mathcal{H}_{V_\varphi, \widetilde{h}}^{\Gamma}(I_0^*)\left(\sum_{i=d_1+1}^{d_2}e^{-\mathbb{Q}S^{\beta_i}}+\sum_{i=d_2+1}^{d}e^{-y\mathbb{Q}T^{\beta_i}}\right)\nonumber\\
&& \ \ =S^{-d_1}\mathcal{H}_{V_\varphi, h_0^*}^{\Gamma}(\widetilde{I}_0)\left(1+\sum_{i=d_1+1}^{d_2}e^{-\mathbb{Q}S^{\beta_i}}\right).
\EQN
Note further  that
\BQN\label{pic1}
\mathcal{H}_{V_\varphi, h^*_0}^{\sup}(\widehat{I}_0)=\mathcal{H}_{V_{\varphi_1}}([0,\sqrt{S}])
\prod_{i=2}^{d_1}\mathcal{H}_{V_{\varphi_i}}[0,S]\prod_{i=d_1+1}^{d_2}
\mathcal{P}_{V_{\varphi_i}}^{h_{i}}[0,S]\rd{\mathcal{H}_{\widetilde{V}_{\varphi, \widetilde{h}}}^{\Gamma^*}(E)}
\EQN
and
\BQN
\mathcal{H}_{V_\varphi, h_0^*}^{\Gamma}(\widetilde{I}_0)=
\prod_{i=1}^{d_1}\mathcal{H}_{V_{\varphi_i}}[0,S]\prod_{i=d_1+1}^{d_2}
\mathcal{P}_{V_{\varphi_i}}^{h_{i}}[0,S]\rd{\mathcal{H}_{\widetilde{V}_{\varphi, \widetilde{h}}}^{\Gamma^*}(E)},
\EQN
with $V_{\varphi_i}, \widetilde{V}_{\varphi}$ and $\widetilde{h}$ defined in \rd{(\ref{vf}) and (\ref{vf1})}. Using further the fact that \cl{(see e.g., Theorem 3.1 in \cite{debicki2002ruin})}$$\lim_{S\rw\IF}\frac{\mathcal{H}_{V_{\varphi_i}}[0,S]}{S}=\mathcal{H}_{V_{\varphi_i}}\in (0,\IF), \quad 1\leq i\leq d_1$$
and letting $S\rw\IF$ on {the left side} of (\ref{main}), we have
\BQNY
\prod_{i=1}^{d_1}\mathcal{H}_{V_{\varphi_i}}\prod_{i=d_1+1}^{d_2}
\lim_{S\rw\IF}\mathcal{P}_{V_{\varphi_i}}^{h_{i}}[-S,S]\mathcal{H}_{\widetilde{V}_{\varphi, \widetilde{h}}}^{\Gamma^*}(E)\leq S^{-d_1}\mathcal{H}_{V_\varphi, h_0^*}^{\Gamma}(\widetilde{I}_0)\left(1+\sum_{i=d_1+1}^{d_2}e^{-\mathbb{Q}S^{\beta_i}}\right)<\IF.
\EQNY
 Thus we conclude that
$$\lim_{S\rw\IF}\mathcal{P}_{V_{\varphi_i}}^{h_{i}}[-S,S]\in (0,\IF), \ \ d_1+1\leq i\leq d_2,$$
which establishes the claim {by letting $S\rw\IF$ on both sides of (\ref{main}).}
For other \HEH{cases of} $a_i, b_i, d_1+1\leq i\leq d_2$, the proof \HEH{is similar as above.}  \QED

	\proofprop{propP} We have that for any $S,T$ positive
	$$ 0<\mathcal{P}_X^{b}(\ehd{[0,S], [0,T]})\leq  \mathcal{P}_X^{b\sigma^2(t)}[0,T].$$
In order to complete the proof it suffices to prove that  $\lim_{T\rw\IF}\mathcal{P}_X^{b\sigma^2(t)}[0,T]<\IF$. 	For this purpose, define for any $S>0, u>1$
$$Y_u(t)=\frac{\overline{X}(u(t+1))}{1+\frac{b\sigma^2(ut)}{2\sigma^2(u)}}, \quad t\in [0, u^{-1}\ln u]. $$
Note that
\BQNY
1-Cor\left(X(ut), X(us)\right)=\frac{\sigma^2(u|t-s|)-(\sigma(ut)-\sigma(us))^2}{2\sigma(ut)\sigma(us)}
=\frac{\sigma^2(u|t-s|)-(u\dot{\sigma}(u\theta)(t-s))^2}{2\sigma(ut)\sigma(us)},
\EQNY
with $\theta\in [s,t]$. By {\bf A1} and Theorem 1.7.2 in \cite{BI1989}, it follows that
$$\lim_{u\rw\IF}\frac{u\dot{\sigma}(u)}{\sigma(u)}=\alpha_\IF.$$
If we set  $f(t)=t^2/\sigma^2(t),$ then by Lemma 5.2 in \cite{KP2015} it follows that   $f$ is  bounded over any compact set and regularly varying at $\IF$ with index $2-2\alpha_\IF>0$. Consequently, UCT implies for any $S>0$
\BQNY
\lim_{u\rw\IF}\sup_{t\in (0, S]}\left|\frac{f(ut)}{f(u)}-|t|^{2-2\alpha_\IF}\right|=0
\EQNY
and therefore as $u\to \IF$
\BQN\label{New}
1-Cor\left(X(ut), X(us)\right)&\sim&\frac{\sigma^2(u|t-s|)}{2\sigma(ut)\sigma(us)}\left(1-\frac{\alpha_\IF^2}{\theta^2}\frac{\sigma^2(u\theta)
(t-s)^2}{\sigma^2(u|t-s|)}\right)\nonumber\\
&=&\frac{\sigma^2(u|t-s|)}{2\sigma(ut)\sigma(us)}\left(1-\alpha_\IF^2\frac{f(u|t-s|)}{f(u\theta)}\right)\sim \frac{\sigma^2(u|t-s|)}{2\sigma^2(u)}
\EQN
for $s,t\in [1,1+u^{-1}\ln u]$. Let further
$$I_k(u)=[ku^{-1}S, u^{-1}(k+1)S ], \quad 0\leq k\leq N(u), \text{  with } N(u):=[S^{-1}\ln u]+1.$$  It follows that for $S$ sufficiently large
\BQN\label{New1}
p_0(u)\leq \pk{\sup_{t\in [0, u^{-1}\ln u]}Y_u(t)>\sqrt{2}\sigma(u)}\leq p_0(u)+\sum_{k=1}^{N(u)}p_k(u),
\EQN
 where
$$p_0(u)=\pk{\sup_{t\in I_0(u)}Y_u(t)>\sqrt{2}\sigma(u)}, $$
$$ p_k(u)=\pk{\sup_{t\in I_k(u)}\overline{X}(u(t+1))>\sqrt{2}\sigma(u)\left(1+\frac{b\sigma^2(kS)}{4\sigma^2(u)}\right)}, \quad  k\geq 1. $$
In order to apply Theorem \ref{lem51}, in view of  (\ref{New}) we set (using the notation in Theorem \ref{lem51})
\BQN K_u=\{k: 0\leq k\leq N(u)\}, \ \ E=[0,S], \ \ g_{u,k}=\sqrt{2}\sigma(u)\left(1+\frac{b\sigma^2(kS)}{4\sigma^2(u)}\right), k\in K_u,
\EQN $$Z_{u,k}(t)=\overline{X}(u(u^{-1}kS+u^{-1}t+1)), \quad k\in K_u,$$
$$\theta_{u,k}(s,t)=g_{u,k}^2\frac{\sigma^2(|t-s|)}{2\sigma^2(u)},\quad  s,t \in E, k\in K_u,$$
$$h_{u,0}(t)=\frac{b\sigma^2(t)}{2\sigma^2(u)},\quad  t\in E, \ \ h_{u,k}=0, \quad k\in K_u \setminus \{0\}, \ \ \eta=X. $$
 {\bf C0} and {\bf C2} are obviously fulfilled.  {\bf C1} is also satisfied with $$ g_{u,0}^2h_{u,0}(t)\rw b\sigma^2(t), \quad u\rw\IF
$$
uniformly with respect to $t\in E$ and
$$g_{u,k}^2h_{u,k}(t)=0, \ \ t\in E, k\in K_u \ehd{\setminus} \{0\},  \quad u>0
$$
 Next  we shall verify  {\bf C3}. Clearly by {\bf A2} for $u$ sufficiently large
$$\theta_{u,k}(s,t)=g_{u,k}^2\frac{\sigma^2(|t-s|)}{2\sigma^2(u)}\leq 2\sigma^2(|t-s|)\leq Q|t-s|^{\alpha_0}, \ \ s,t\in E, k\in K_u.$$
Moreover, by (\ref{New})
\BQNY
&&\sup_{k \in K_u}\sup_{\norm{t-s}<\epsilon,s,t\in E}g_{u,k}^2
\EE{ \left[ Z_{u,k}(t)-Z_{u,\tau }(s)\right]Z_{u,k}(0)}\\
&& \ \ \leq \sup_{k \in K_u}\sup_{\norm{t-s}<\epsilon,s,t\in E}g_{u,k}^2
\left(\frac{\sigma^2(t)}{2\sigma^2(u)}(1+o(1))-\frac{\sigma^2(s)}{2\sigma^2(u)}(1+o(1))\right)\\
&& \ \ \leq \sup_{k \in K_u}\sup_{\norm{t-s}<\epsilon,s,t\in E}\frac{g_{u,k}^2}{2\sigma^2(u)}
\left(\abs{\sigma^2(t)-\sigma^2(s)}+o(1)\right)\rw 0, \ \  u\rw\IF, \epsilon\downarrow 0.
\EQNY
Thus {\bf C3} is satisfied. Therefore, in light of Theorem \ref{lem51}, we have that
\BQNY
\limit{u} \frac{p_0(u)}{\Psi(\sqrt{2}\sigma(u))}=  \mathcal{P}_X^{b\sigma^2(t)}[0,S]
\EQNY
and
\BQN
\lim_{u\rw\IF}\sup_{k\in  K_u/\{0\}}\ABs{\frac{p_k(u)}{\Psi\left(\sqrt{2}\sigma(u)\left(1+\frac{b\sigma^2(kS)}{4\sigma^2(u)}\right)\right)}-\mathcal{H}_X[0,S]}=0.
\EQN
Dividing (\ref{New1}) by $\Psi(\sqrt{2}\sigma(u))$, letting $u\rw\IF$  and by {\bf A1}, we have that  for  sufficiently large $S_1$
\BQNY
\mathcal{P}_X^{b\sigma^2(t)}[0,S]&\leq& \mathcal{P}_X^{b\sigma^2(t)}[0,S_1]+\mathcal{H}_X[0,S_1]\sum_{k=1}^{\IF} e^{-\frac{b\sigma^2(kS_1)}{2}}\\
&\leq& \mathcal{P}_X^{b\sigma^2(t)}[0,S_1]+\mathcal{H}_X[0,S_1]\sum_{k=1}^{\IF} e^{-Q_1(kS_1)^{\alpha_\IF}}\\
&\leq & \mathcal{P}_X^{b\sigma^2(t)}[0,S_1]+\mathcal{H}_X[0,S_1] e^{-Q_2S_1^{\alpha_\IF}}.
\EQNY
Next, letting $S\rw\IF$ leads to
\BQNY
\lim_{S\rw\IF}\mathcal{P}_X^{b\sigma^2(t)}[0,S]\leq \mathcal{P}_X^{b\sigma^2(t)}[0,S_1]+\mathcal{H}_X[0,S_1] e^{-Q_2S_1^{\alpha_\IF}}<\IF
\EQNY
establishing the claim.

\QED

\section{Appendix}
\proofrem{remark} ii). First we suppose that {\bf C2} and (\ref{eq2}) hold. Our aim is to prove (\ref{eq3}).
By (\ref{eq2}), the continuity of $\sigma^2_\eta(t), t\in\bD$ and the compactness of $\bD$, for any $c>0$, there exists a constant $\gE{\epsilon}:= \epsilon_c>0$ such that
$$\limsup_{u\rw\IF}\sup_{\gE{\tau_u}\in K_u}\sup_{\norm{t-s}< \gE{\epsilon},s,t\in E}\ABs{ \gtu^2Var\left( \gE{\qqq_u(s)}\right)-\gtu^2Var\left(\gE{\qqq_u(t)}\right)}<c/3,$$
with $\gE{  \qqq_u(t)= Z_{u,\tau_u}(t)- Z_{u,\tau_u}(0)}$ and further
$$\sup_{\norm{t-s}<\gE{\epsilon},s,t\in E}\ABs{ \sigma_\eta^2(t)-\sigma^2_\eta(s)}<c/3.$$
By the compactness of $\bD$, we can find $\bD_c\subset E$ which has a finite number of elements  such that for any $t\in\bD$
$$O_{\gE{\epsilon}}(t)\cap \bD_c\neq \emptyset, \ \ \ O_{\gE{\epsilon}}(t):=\{s\in \mathbb{R}^d : \norm{t-s}<
\gE{\epsilon}\}.$$
 For any $t\in \bD$, with $t'\in O_{\gE{\epsilon}}(t)\cap \bD_c$
\BQNY
\ABs{ \gtu^2Var\left(\gE{  \qqq_u(t)}\right)-2\sigma_\eta^2(t)}&\leq& \ABs{ \gtu^2Var\left(\gE{  \qqq_u(t)}\right)-\gtu^2Var\left(\gE{  \qqq_u(t')}\right)}\\
&&+2\ABs{ \sigma_\eta^2(t)-\sigma^2_\eta(t')}
+\ABs{ \gtu^2Var\left(\gE{  \qqq_u(t')}\right)-2\sigma_\eta^2(t')}.
\EQNY
It follows from {\bf C2} that
$$\lim_{u\rw\IF} \sup_{ \gE{\tau_u}\in K_u} \ABs{ g_{u,\tau_u}^2Var(\gE{  \qqq_u(t)})- 2 \sigma_\eta^2(t)}=0, \ \ t\in E.$$
Consequently, we have
\BQNY
&&\limsup_{u\rw\IF} \sup_{\gE{\tau_u}\in K_u}\sup_{ t\in \bD} \ABs{ \gtu^2Var\left(\gE{  \qqq_u(t)}\right)-2\sigma_\eta^2(t)}\\
&& \ \  \leq
\limsup_{u\rw\IF}\sup_{\gE{\tau_u}\in K_u}\sup_{\norm{t-s}<\gE{\epsilon},s,t\in E}\ABs{ \gtu^2Var\left(
\gE{  \qqq_u(s)}\right) -\gtu^2Var\left(\gE{  \qqq_u(t)}\right)}\\
&&  \ \ \ +2\sup_{\norm{t-s}< \gE{\epsilon},s,t\in E}\ABs{ \sigma_\eta^2(t)-\sigma^2_\eta(s)}
 +\limsup_{u\rw\IF}\sup_{\gE{\tau_u}\in K_u}\sup_{t\in \bD_{\gE{c}}}\ABs{ \gtu^2Var\left(\gE{  \qqq_u(t)}\right)-2\sigma_\eta^2(t)}\\
&& \ \ \ \le  c.
\EQNY
{\kd{Hence} letting $c$ to 0 yields} (\ref{eq3}). \\
Next,  supposing that {\bf C2} and (\ref{eq3}) hold, we prove (\ref{eq2}).
By  the continuity of $\sigma^2_\eta(t), t\in\bD$ and the compactness of $\bD$, for any $c>0$, there exists a constant $\gE{\epsilon}>0$ such that
$$\sup_{\norm{t-s}<\gE{\epsilon},s,t\in E}\ABs{ \sigma_\eta^2(t)-\sigma^2_\eta(s)}<c/3.$$
{For any $s,t \in E$}
\BQNY
\ABs{ \gtu^2Var\left(\gE{  \qqq_u(s)}\right)-\gtu^2Var\left(\gE{  \qqq_u(t)}\right)}
&\leq & \ABs{ \gtu^2Var\left(\gE{  \qqq_u(s)}\right)-2\sigma^2_\eta(s)}+2|\sigma^2_\eta(s)-\sigma^2_\eta(t)|\\
&& +\ABs{2\sigma_\eta^2(t)-\gtu^2Var\left(\gE{  \qqq_u(t)}\right)}.
\EQNY
Consequently, by (\ref{eq3})
\BQNY
\lefteqn{\limsup_{u\rw\IF}\sup_{\gE{\tau_u}\in K_u}\sup_{\norm{t-s}<\gE{\epsilon},s,t\in E}\ABs{ \gtu^2Var\left(\gE{  \qqq_u(s)}\right)-\gtu^2Var\left(\gE{  \qqq_u(t)}\right)}}\\
 &  \leq &  2\limsup_{u\rw\IF}\sup_{\gE{\tau_u}\in K_u}\sup_{t\in E}\ABs{ \gtu^2Var\left(\gE{  \qqq_u(t)}\right)-2\sigma_\eta^2(t)}+2\sup_{\norm{t-s}<\gE{\epsilon},s,t\in E}\ABs{ \sigma_\eta^2(t)-\sigma^2_\eta(s)}\\
 & \leq &  c.
\EQNY
\kd{Letting $c\to0$, the above establishes
\eqref{eq2}, which completes the proof.} \QED

{\bf Acknowledgement}: We would like to thank the referees
for their useful \kk{comments} 
leading to significant improvement for the readability of this paper.  Thanks to Swiss National Science Foundation grant No. 200021-166274.  KD acknowledges partial support by NCN Grant No 2015/17/B/ST1/01102 (2016-2019).

\bibliographystyle{plain}

\bibliography{UnifJAP13}

\end{document}